\documentclass[12pt]{amsart}
\usepackage{latexsym}
\usepackage{amssymb}
\usepackage{amscd}
\renewcommand\a{\alpha}
\renewcommand\b{\beta}
\newcommand\g{\gamma}
\renewcommand\d{\delta}
\newcommand\la{\lambda}
\newcommand\z{\zeta}
\newcommand\e{\eta}
\renewcommand\th{\theta}
\newcommand\io{\iota}

\newcommand\n{\nu}
\newcommand\s{\sigma}
\newcommand\x{\chi}
\newcommand\f{\phi}
\newcommand\vf{\varphi}
\newcommand\p{\psi}

\renewcommand\r{\rho}

\newcommand\vS{\varSigma}

\newcommand\F{\Phi}

\newcommand\vG{\varGamma}
\newcommand\ve{\varepsilon}

\newcommand{\ZZ}{\mathbb Z}

\newcommand{\HH}{\mathbb H}
\renewcommand{\AA}{\mathbb A}

\newcommand\Fq{{\mathbf F}_q}

\newcommand\Ql{\bar{\mathbf Q}_l}
\newcommand\BA{\mathbf A}
\newcommand\BP{\mathbf P}

\newcommand\BF{\mathbf F}
\newcommand\BC{\mathbf C}

\newcommand\BS{\mathbf S}
\newcommand\BZ{\mathbf Z}

\newcommand\BH{\mathbf H} 

\newcommand\BG{\mathbf G}

\newcommand\Br{\mathbf r}

\newcommand\CB{\mathcal{B}}
\newcommand\ZC{\mathcal{C}}
\newcommand\CH{\mathcal{H}}

\newcommand\CE{\mathcal{E}}
\newcommand\CL{\mathcal{L}}

\newcommand\CS{\mathcal{S}}
\newcommand\CM{\mathcal{M}}
\newcommand\CN{\mathcal{N}}

\newcommand\CP{\mathcal{P}}
\newcommand\CV{\mathcal{V}}
\newcommand\CF{\mathcal{F}}
\newcommand\CG{\mathcal{G}}

\newcommand\CW{ \mathcal{W}}

\newcommand\Fg{\mathfrak g}

\newcommand\Fh{\mathfrak h}
\newcommand\Fs{\mathfrak s}
\newcommand\Fl{\mathfrak l}

\newcommand\Fm{\mathfrak m}
\newcommand\Fn{\mathfrak n}
\newcommand\Fp{\mathfrak p}
\newcommand\Fz{\mathfrak z}
\newcommand\iv{^{-1}}
\newcommand\wh{\widehat}
\newcommand\wt{\widetilde}
\newcommand\wg{^{\wedge}}

\newcommand\ol{\overline}
\newcommand\hra{\hookrightarrow}
\newcommand\lra{\leftrightarrow}

\newcommand\bs{\backslash}
\newcommand\IC{\operatorname{IC}}
\newcommand\Ker{\operatorname{Ker}}
\newcommand\Hom{\operatorname{Hom}}
\newcommand\End{\operatorname{End}}
\newcommand\Aut{\operatorname{Aut}}

\newcommand\Lie{\operatorname{Lie}}
\newcommand\Tr{\operatorname{Tr}\,}

\newcommand\pr{\operatorname{pr}} 
\newcommand\Ad{\operatorname{Ad}}
\newcommand\ad{\operatorname{ad}}

\newcommand\reg{_{\operatorname{reg}}}
\newcommand\uni{_{\operatorname{uni}}}
\newcommand\nil{_{\operatorname{nil}}}
\newcommand\id{\operatorname{id}}

\newcommand\lp{\operatorname{\!\langle\!}}
\newcommand\rp{\operatorname{\!\rangle\!}}
\renewcommand\Im{\operatorname{Im}}

\newcommand\odd{\operatorname{odd}}

\newcommand\dw{\dot w}

\newcommand{\isom}{\,\raise2pt\hbox{$\underrightarrow{\sim}$}\,}
\numberwithin{equation}{section}

\newtheorem{thm}{Theorem}[section]
\newtheorem{lem}[thm]{Lemma}
\newtheorem{cor}[thm]{Corollary}
\newtheorem{prop}[thm]{Proposition}

\def \para#1{\par\medskip\textbf{#1}
              \addtocounter{thm}{1}}
\begin{document}
\setlength{\baselineskip}{4.9mm}
\setlength{\abovedisplayskip}{4.5mm}
\setlength{\belowdisplayskip}{4.5mm}
\renewcommand{\theenumi}{\roman{enumi}}
\renewcommand{\labelenumi}{(\theenumi)}
\renewcommand{\thefootnote}{\fnsymbol{footnote}}
\renewcommand{\thefootnote}{\fnsymbol{footnote}}
\parindent=20pt
\medskip
\title{Generalized Green functions \\and unipotent classes
\\ for finite reductive groups, I }
\author{Toshiaki Shoji}

\maketitle
\begin{center}
Graduate School of Mathematics \\
Nagoya University  \\
Chikusa-ku, Nagoya, 464-8602,  Japan
\end{center}
\pagestyle{myheadings}
\markboth{SHOJI}{GENERALIZED GREEN FUNCTIONS}

\begin{abstract}
The algorithm of computing generalized Green functions
of a reductive group $G$ contains some unknown scalars
occurring from the $\Fq$-structure of irreducible local 
systems on unipotent classes of $G$.  
In this paper, we determine such scalars in the case where
$G = SL_n$ with Frobenius map $F$ of split type 
or non-split type.  In the case where $F$ is of non-split type,
we use the theory of graded Hecke algebras due to Lusztig. 
\end{abstract}

\bigskip
\medskip
\addtocounter{section}{-1}
\section{Introduction}
Let $G$ be a connected reductive group defined over a finite 
field $\Fq$ with Frobenius map $F$.  
In [L1], Lusztig classified the irreducible characters of
finite reductive groups $G^F$ in the case where the center of 
$G$ is connected.  Later in [L5], he extended his results to the 
disconnected center case.  In the course of the classification, in 
particular in the connected center case, he defined almost characters
of $G^F$, which forms an orthonormal basis of the space $\CV(G^F)$ of class
functions of $G^F$ different from the basis consisting of 
irreducible characters.  They are
defined as explicit linear combinations of irreducible characters,
and the transition matrix between these two bases are almost diagonal.
So, the determination of the character values of irreducible
characters of $G^F$ is equivalent to that of almost
characters.
\par
On the other hand, Luszitg founded in [L3] the theory of character 
sheaves, and showed that the characteristic functions of character
sheaves form an orthonormal basis of $\CV(G^F)$.  He conjectured 
that those functions coincide, 
up to scalar, with almost characters (with an appropriate 
generalization of almost characters if the center is disconnected).  
Lusztig's conjecture was proved by 
the author in [S3] in the case where the center is connected.  It was 
also proved for certain groups with disconnected center, i.e., 
for $Sp_{2n}$ and (under a suitable modification for a disconnected
group) $O_{2n}$ with ch $\Fq \ne 2$ by 
Waldspurger [W],   
for $SL_n$ by the author [S4] (with ch $\Fq$ not too samll), 
and independently, for $SL_n$ and $SU_n$ by Bonnaf\'e [B] 
(with $q$ not too small.) 
\par
If Lusztig's conjecture is established, the computation of irreducible
characters of $G^F$ is reduced to the computation of characteristic 
functions of character sheaves, and to the determination of 
scalars involved in Lusztig's conjecture.  
In [L3], Lusztig proved that the computation of
the characteristic functions of character sheaves are reduced to
the computation of generalized Green functions of various reductive
subgroups of $G^F$.  Then he showed that there exists a general 
algorithm of computing generalized Green functions. More precisely, 
he showed that generalized Green functions can be expressed as an 
explicit linear combination of various characteristic functions 
$\x_{C',\CE'}$ of
the $G$-equivariant local system $\CE'$ on a unipotent class $C'$ 
in $G$.  Up to scalar, $\x_{C',\CE'}$ can be easily described in terms
of the irreducible character of the component group 
$A_G(u) = Z_G(u)/Z_G^0(u)$ for $u \in {C'}^F$ corresponding to $\CE'$. 
However, this scalar depends on the choice of the isomorphism 
$F^*\CE \isom \CE$ for a cuspidal pair $(C,\CE)$ on a Levi subgroup
$L$ of a parabolic subgroup $P$ of $G$, and on the intersection 
cohomology complex $K$ induced from $\CE\boxtimes\Ql$ on 
$C\times Z_L^0$ (see (1.2.2)). 
\par
The purpose of this paper is to determine these scalars occurring 
in the computation of generalized Green functions.  In the case 
of Green functions, this problem is equivalent to determining 
a representative $u \in {C'}^F$ such that the action of 
$F$ on the $l$-adic cohomology group $H^{m}(\CB_u, \Ql)$ can 
be described explicitly, where $\CB_u$ is the variety of Borel
subgroups of $G$ containing $u$, and $m/2 = \dim \CB_u$.
It was shown in [S1], [S2] and [BS] that there exists a unipotent 
element $u \in {C'}^F$, in the case where $G^F$ is of split type, and 
$G$ is not of type $E_8$, such that $F$ acts on $H^{m}(\CB_u,\Ql)$
by a scalar multiplication $q^{m/2}$.  Such a unipotent element is 
called a split element.  Even if the remaining cases, the action of 
$F$ can be described, and by using this, Green functions of exceptional
groups ($F_4, E_6, E_7$ and $E_8$) were computed explicitly 
by [S1], [BS] for a good characteristic case.  The case $G_2$ had been
computed by Springer [Spr] in an earlier stage.  
(Green functions of exceptional groups in certain bad characteristic 
case were computed by Malle [M] by a direct computation).
\par
In the case of generalized Green functions, one has to consider the
cohomology group $H^{m}_c(\CP_u, \dot\CE)$, where $\CP_u$ is a
certain subvariety of parabolic subgroups of $G$ conjugate to $P$, 
and $\dot\CE$ is a local system on $\CP_u$ determined from the cuspidal pair 
$(C,\CE)$ on a Levi subgroup $L$ of $P$, and $m/2 = \dim \CP_u$.     
We need to describe the action of $F$ on such cohomology groups.
This problem is reduced to the case where $G$ is simply connected, and 
simple modulo center. In this paper, we discuss the case where 
$G = SL_n$ with $F$ of split type or non-split type. 
In the case where $F$ is of split type, the method employed here
is to compare the Frobenius action in the case of 
$SL_n$ with $SL_{n-1}$, 
which is a natural generalization of the method in the case of $GL_n$.
In the case of $GL_n$ with $F$ of non-split type, the Frobenius action
was determined by investigating the action of $F$ on $H^*(\CB, \Ql)$ 
by making use of the $F$-equivariant surjective map 
$\pi_u: H^{m}(\CB,\Ql) \to H^{m}(\CB_u, \Ql)$
induced from the inclusion $\CB_u \hra \CB$, where $\CB$ is the 
Flag variety of $G$.  
However, this argument is not generalized to our case.
Although we have a counter part $\CP_{u_1}$ of $\CB$, and a natural
map $\pi_u: H_c^{m}(\CP_{u_1},\dot\CE) \to H^{m}_c(\CP_u,\dot\CE)$, 
there does not exist an 
immersion $\CP_u \hra \CP_{u_1}$, and the surjectivity of $\pi_u$ 
is no longer trivial.  In order to overcome such difficulties, following 
the idea of Lusztig, we appeal to the theory of graded Hecke algebra
developed in [L7], which makes it possible to compare the Frobenius
actions via the isomorphism  
$H^0_c(\CP_{u_1},\dot\CE) \simeq H_c^0(\CP_u,\dot\CE)\simeq \Ql$. 
\par
The remaining cases where $G \ne SL_n$ will be treated in a 
subsequent paper.
\par
The author is grateful to G. Lusztig for stimulating discussions on 
graded Hecke algebras.  

\section{Preliminaries}
\para{1.1.}
Let $G$ be a connected reductive algebraic group 
over a field $k$, where $k$ is an algebraic closure of a finite
filed $\Fq$ of characteristic $p$.
Let $C$ be a unipotent conjugacy class in $G$, and $\CE$  
an irreducible local system on $C$ which is $G$-equivariant for
the conjugation action.
$\CE$ is called a cuspidal local system on $C$ if the following
condition is satisfied: for any proper parabolic subgroup $P$ 
of $G$ with Levi decomposition $P = LU_P$ and for any unipotent element
$u \in L$, we have $H_c^{\d}(uU_P \cap C, \CE) = 0$,
where $\d = \dim C - \dim \text{(class of $u$ in L)}$ (cf. [L2, 2.4]).
It is known by [L3, V, 23.1], that if $p$ is almost good then 
the above condition is
equivalent to the condition that $H_c^i(uU_P \cap C, \CE) = 0$ for
any $i$ (i.e., $\CE$ is strongly cuspidal).  We also say that 
$(C, \CE)$ is a cuspidal pair in $G$.
\par
Let $\CN_G$ be the set of pairs $(C', \CE')$ up to $G$-conjugacy, 
where $C'$ is a unipotent class in $G$ and $\CE'$ is a
$G$-equivariant irreducible local system on $C$.  We also denote
by $\CM_G$ the set of triples $(L, C, \CE)$ up to $G$-conjugacy,
where $L$ is a Levi subgroup of a parabolic subgroup of $G$, and
$\CE$ is a cuspidal local system on a unipotent class $C$ of $L$. 
In [L2, 6.5], Lusztig  has shown that 
there exists a natural bijection 
\begin{equation*}\tag{1.1.1}
\CN_G \simeq \coprod_{(L,C,\CE) \in \CM_G} (N_G(L)/L)\wg, 
\end{equation*}
which is called the generalized Springer correspondence between
unipotent classes and irreducible characters of various Coxeter groups.
(For a finite group $H$, we denote by $H\wg$ the set of irreducible
characters of $H$). 
Note that $N_G(L)/L$ is a Coxeter group with standard generators 
whenever $(L,C,\CE) \in \CM_G$.
\para{1.2.}
We describe the generalized Springer correspondence more precisely.
Take $(L,C,\CE) \in \CM_G$.  
Let $Z^0_L$ be the connected center of $L$, and put 
$\wt C\reg = C\cdot(Z^0_L)\reg \subset  \wt C = C\cdot Z^0_L $, where
\begin{equation*}
(Z^0_L)\reg = \{ z \in Z_L^0 \mid Z^0_G(z) = L \}.
\end{equation*}
We define a diagram
\begin{equation*}
\begin{CD}
\wt C @<\a_1<< \wh Y @>\b_1>> \wt Y @>\pi>> Y,
\end{CD}
\tag{1.2.1}
\end{equation*}
where 
\begin{align*}
Y &= \bigcup_{x \in G}x\wt C\reg x\iv \subset G, \\
\wt Y &= \{ (g, xL) \in G \times (G/L) \mid x\iv gx \in \wt C\reg \}, \\
\wh Y &= \{ (g,x) \in G \times G \mid x\iv gx \in \wt C \},
\end{align*}
and
\begin{equation*}
\a_1(g,x) = x\iv gx, \quad \b_1(g,x) = (g,xL), \quad \pi(g,xL) = g.
\end{equation*}
Then $Y$ is a smooth, irreducible subvariety of $G$, and $\pi$
is a principal covering of $Y$ with group $\CW = N_G(L)/L$.
There is a canonical local system $\wt\CE$ on $\wt Y$ satisfying the
property that $\b_1^*\wt\CE = \a_1^*(\CE\boxtimes \Ql)$, where 
$\CE\boxtimes \Ql$ is the inverse image of $\CE$ under the natural 
map $\wt C  = C\times Z^0_L \to C$.  We define an intersection 
cohomology complex $K$ by
\begin{equation*}
\tag{1.2.2}
K = \IC(\ol Y, \pi_!\wt\CE)[\dim Y]
\end{equation*}
and regard it as a perverse sheaf on $G$ by extending by 0 outside of
$\ol Y$. 
Lusztig showed that $K$ is a $G$-equivariant semisimple perverse sheaf on $G$, 
and that $\End K \simeq \Ql[\CW]$.
It follows that $K$ can be decomposed as 
\begin{equation*}
\tag{1.2.3}
K \simeq \bigoplus_{E \in \CW\wg}V_E\otimes K_E,
\end{equation*}
where $K_E$ is a simple perverse sheaf on $G$ such that  
$V_E = \Hom (K_E, K)$ is 
an irreducible $\CW$-module corresponding to $E \in \CW\wg$.  
\par
Let $G\uni$ be the unipotent variety of $G$.
Then $K[-d]|_{G\uni}$ turns out to be a $G$-equivariant semisimple 
 perverse sheaf on $G\uni$, where 
$d = \dim Z^0_L = \dim Y - \dim (Y \cap G\uni)$.  
Hence it is decomposed as 
\begin{equation*}
\tag{1.2.4}
K[-d]|_{G\uni} = \bigoplus_{(C',\CE') \in \CN_G}
                     V_{(C',\CE')}\otimes \IC(\ol C', \CE')[\dim C'],
\end{equation*}
where $V_{(C',\CE')}$ is a multiplicity space for the simple 
perverse sheaf $\IC(\ol C',\CE')[\dim C']$ on $G\uni$.
Comparing (1.2.3) with (1.2.4), we see that for each $E \in \CW\wg$, 
there exists a pair $(C',\CE') \in \CN_G$ such that
\begin{equation*}
\tag{1.2.5}
K_E|_{G\uni} \simeq \IC(\ol C',\CE')[\dim C' + \dim Z^0_L].
\end{equation*}
The correspondence $E \mapsto (C', \CE')$ gives a bijection
$\coprod_{(L,C,\CE)} (N_G(L)/L)\wg \to \CN_G$ in (1.1.1).
\para{1.3.}
We now consider the $\Fq$-structure on $G$. So assume that 
$G$ is defined over $\Fq$ with Frobenius endomorphism 
$F: G \to G$.
Then $F$ acts naturally on the set $\CN_G$ and $\CM_G$ by 
$(C',\CE') \mapsto (F\iv C', F^*\CE')$, 
  $(L, C, \CE) \mapsto (F\iv L, F\iv C, F^*\CE)$, and the map
in (1.1.1) is compatible with $F$-action.
Now assume that $(L,C,\CE) \in \CM_G$ is $F$-stable.
Then we may choose $(L,C, \CE)$, as a representative of its 
$G$-conjugacy class, such that $L$ is an $F$-stable Levi subgroup
of an $F$-stable parabolic subgroup $P$ of $G$, with 
$FC = C, F^*\CE \simeq \CE$. 
We choose an isomorphism $\vf_0: F^*\CE \isom \CE$
which induces a map of finite order on the stalk of $\CE$ at any
point of $C^F$.  Since the diagram in (1.2.1), and so the construction 
of the complex $K$ is compatible with 
$\Fq$-structure, $\vf_0$ induces a natural isomorphism 
$\vf: F^*K \isom K$.  
We consider the characteristic function $\x_{K,\vf}$ of $K$.
The restriction of $\x_{K,\vf}$ on $G\uni$ gives a $G^F$-invariant
function on $G^F\uni$, which is the generalized
Green function $Q^G_{L,C,\CE,\vf_0}$ (cf. [L3, II]). 
\par
Here $F$ acts naturally on $\CW$, which induces a Coxeter group
automorphism of degree, say $c$.  We consider the 
semidirect product $\wt\CW = \CW \rtimes (\BZ/c\BZ)$. 
If an irreducible representation $V_E$ of $\CW$ is $F$-stable, 
it can be extended to an irreducible representation of $\wt\CW$, 
in $c$ different ways.
Assume that $E \in \CW\wg$ is $F$-stable.  Then the corresponding 
$(C',\CE') \in \CN_G$ is
also $F$-stable, and we have $F^*K_E \isom K_E$.
A choice of an isomorphism $\vf_E: F^*K_E \isom K_E$ induces a
bijection $\s_E: V_E \to V_E$, which makes $V_E$ into an irreducible
$\wt\CW$-module $\wt V_E$.  We choose $\vf_E$ so that $\wt V_E$ turns
out to be a preferred extension of $V_E$ (cf. [L3, IV, (17.2)].
By making use of $\vf_E: F^*K_E \simeq K_E$, we shall define an
isomorphism $\p: F^*\CE' \isom \CE'$ as follows; By (1.2.5), we have
$\CH^{a_0}(K_E)|_{C'} = \CE'$ for $a_0 = -\dim Z_L^0 - \dim C'$.
We define $\p$ so that $q^{(a_0+r)/2}\p$ corresponds to the map
defined by $\vf_E: F^*\CH^{a_0}(K_E) \isom \CH^{a_0}(K_E)$, where 
\begin{equation*}
r = \dim Y = \dim G - \dim L + \dim (C \times Z_L^0),
\end{equation*}
and so
\begin{equation*}
\tag{1.3.1}
a_0+r = (\dim G - \dim C') - (\dim L - \dim C).
\end{equation*}
\par
We define a function $Y_j$ on $G^F\uni$ 
for each $j = (C',\CE') \in \CN_G^F$
by
\begin{equation*}
Y_j(g) = \begin{cases}
                     \Tr(\p, g) &\quad\text{ if } g \in {C'}^F, \\
                      0         &\quad\text{ if } g \notin {C'}^F.
                   \end{cases}
\end{equation*}
Then $\{ Y_j \mid j \in \CN_G^F\}$ gives rise to a basis of the 
space of $G^F$-invariant functions on $G^F\uni$.
Now the computation of $\x_{K,\vf}$ is reduced to the computation of
$\x_{K_E, \vf_E}$ for each $F$-stable irreducible character $\CE$ of $\CW$.
We denote $\x_{K_E,\vf_E}$ by $X_j$ if $E$ corresponds to
$j = (C',\CE')$ under the generalized Springer correspondence.
In [L3, V], Luszitg gave a general algorithm of expressing $X_i$
as an explicit linear combination of various $Y_j$. 
Thus the computation of $\x_{K, \vf}$ is reduced to the computation of
$Y_j$.
\par
We shall describe the functions $Y_j$. Let us choose 
$u \in {C'}^F$, and put $A_G(u) = Z_G(u)/Z_G^0(u)$.  Then $F$
acts naturally on $A_G(u)$, and the set of $G$-equivariant simple local 
systems on $C'$ is in bijective correspondence with the set of 
$F$-stable irreducible characters of $A_G(u)$. 
Let us denote by $\r$ the irreducible character of $A_G(u)$ 
corresponding to $\CE'$.  Let $\s$ be the restriction of $F$ on
$A_G(u)$.  Then $\r$ can be extended to an irreducible character
of the semidirect group $\wt A_G(u) = A_G(u)\rtimes \lp\s\rp$.  
We choose an extension $\wt\r$ of $\r$.
$\CE'_u$ has a structure of $A_G(u)$-module affording the
character $\r$, which is extended to the $\wt A_G(u)$-module affording
$\wt\r$.  We choose an isomorphism 
$\p_0: F^*\CE' \isom \CE'$ by the condition that $\p_0$ induces
an isomorphsim on $\CE'_u$ corresponding to the action of $\s$ on 
$\wt\r$.  
\par
Since $\CE'$ is a simple local system, there exists $\g \in \Ql^*$
(depending on the choice of $\vf_0$, $u$ and $\wt\r$) such that 
$\p = \g\p_0$.  We define functions $Y_j^0$ on the set $G^F\uni$ 
in a similar way as $Y_j$, but replacing $\p$ by $\p_0$.  Then 
clearly we have $Y_j = \g Y_j^0$. 
We note that the functions $Y_j^0$ are described in an explicit 
way as follows. 
The set of $G^F$-conjugacy classes in
${C'}^F$ is in bijective correspondence with the set of $F$-twisted 
conjugacy classes in $A_G(u)$.  We denote by $u_a$ a representative in
the $G^F$-conjugacy class contained in ${C'}^F$ corresponding to 
an $F$-twisted conjugacy class in $A_G(u)$ containing $a$.  
Then we have
\begin{equation*}
Y_j^0(g) = \begin{cases}
     \wt\r(a\s) &\quad\text{ if $g$ is $G^F$-conjugate to $u_a$,} \\
      0    &\quad\text{ if $g \notin {C'}^F$.}
            \end{cases}
\end{equation*}
\par
It follows from the above discussion that the computation of 
generalized Green functions is reduced to
the determination of the scalar constant $\g$ for each pair 
$(C',\CE') \in \CN_G^F$.  Let us choose $v \in C^F$, and let
$\r_0$ be the $F$-stable irreducible character of $A_L(v)$ 
corresponding to $\CE$.  Then as in the discussion above, the
isomorpshims $\vf_0: F^*\CE \isom \CE$ is given by choosing 
an extension $\wt\r_0$ of $\r_0$ to the semidirect group 
$\wt A_L(v) = A_L(v)\rtimes \lp\s\rp$.  Thus $\g$ is determined by 
$v, \wt\r_0, u, \wt\r$, which we denote by 
$\g = \g(v,\wt\r_0, u,\wt\r)$. 
The purpose of this paper is to describe
the constants $\g(v,\wt\r_0,u,\wt\r)$ explicitly.
\para{1.4.} 
In order to make the Frobenius action more explicit, we shall consider
the following varieties. Put
\begin{align*}
\tag{1.4.1}
\CP_{u} &= \{ gP \in G/P \mid g\iv ug \in CU_P\}, \\
\wh \CP_{u} &= \{ g \in G \mid g\iv ug \in CU_P\},
\end{align*}
and consider the diagram
\begin{equation*}
\begin{CD}
C @<\a << \wh \CP_{u} @> \b>> \CP_{u}
\end{CD}
\tag{1.4.2}
\end{equation*}
with
\begin{equation*}
\a: g \mapsto C\text{-component of }g\iv ug \in CU_P, \quad
\b: g \mapsto gP.
\end{equation*} 
We define a local system $\dot\CE$ on $\CP_u$ by the property that
$\a^*\CE = \b^*\dot\CE$.
Then it is known by [L3, 24.2.5] that 
\begin{equation*}
\tag{1.4.3}
\CH^{a_0}_u(K) \simeq H_c^{a_0+r}(\CP_u, \dot\CE).
\end{equation*}
It is also known by [L2, 1.2 (b)] that $\dim \CP_u \le (a_0+r)/2$.  Since 
the left hand side of (1.4.3) is non-zero by (1.2.5), we see that 
\begin{equation*}
\tag{1.4.4}
\dim \CP_u = (a_0 + r)/2.
\end{equation*}
\par
Since $P$ is $F$-stable, $\CP_u, \wh \CP_u$ are $F$-stable, and the 
diagram in (1.4.2) is compatible with Frobenius maps.  Moreover,
the isomorphism $\vf_0$ induces an isomorphism 
$\dot\vf_0: F^*\dot\CE \isom \dot\CE$.  This induces a linear map
$\F$ on $V = H_c^{a_0+r}(\CP_u, \dot\CE)$. 
By (1.4.3), $\CW$ acts on $V$.  Also
$Z_G(u)$ acts naturally on $V$, where 
$Z^0_G(u)$ acts trivially on it.  Then it induces an action of 
$A_G(u)$, which commutes with the action 
of $\CW$.  Let $\r$ be an $F$-stable irreducible 
character of $A_G(u)$ corresponding to $\CE'$ as in 1.3, and 
$V_{\r}$ the $\r$-isotipic part of $V$.  
Then $\F$ leaves $V_{\r}$ stable. 
The previous discussion 
shows that $V_{\r}$ can be identified with $\wt V_E\otimes\CE_u'$, 
and $\F|_{V_{\r}}$ coincides with $\s_E\otimes q^{(q_0+r)/2}\p$.  
Thus the map $\p$ can be described by investigating $\F$ on 
$H_c^{a_0+r}(\CP_u,\dot\CE)_{\r}$. 
\para{1.5.}
We show that the description of the mixed structure 
$\p: F^*\CE' \to \CE'$ on $C'$ is reduced to the case where
$G$ is simply connected, almost simple.  In fact, let 
$\pi: G \to G' = G/Z^0_G$ be the natural homomorphism.  
Then $\pi$ induces a bijection between $\CM_G$ (resp. $\CN_G$) 
and $\CM_{G'}$ (resp. $\CN_{G'}$) which commutes with their
$\Fq$-structures.  Hence we may assume that $G$ is semisimple.
Let $\wt\pi : \wt G \to G$ be the simply connected covering of $G$.
Then $(L, C, \CE) \mapsto (\wt\pi\iv(L), C, \wt\pi^*\CE)$ gives
a bijection between the set $\CM_G$ and the subset of $\CM_{\wt G}$
on which $\ker\wt\pi$ acts trivially.  Hence the mixed structure
$\vf_0:F^*\wt\pi^*\CE \to \wt\pi^*\CE$ for the pair $(C, \wt\pi^*\CE)$
on $\wt G$ determines the mixed structure for the pair 
$(C, \CE)$ on $G$. Similarly, $\wt\pi$ induces a bijection between the
set $\CN_G$ and the subset of $\CN_{\wt G}$ on which $\ker\wt\pi$ acts
trivially, and so the mixed structure of the pair $(C',\CE')$ on $G$
is determined by the mixed structure of the pair $(C',\wt\pi^*\CE')$
on $\wt G$.  The procedure of determining the mixed structure of 
$(C', \CE')$ from that of $(C,\CE)$ is parallel for $\wt G$ and $G$. 
\par
It follows from the above discussion that we may assume 
$G$ is simply connected,
semisimple.  Then $G$ is isomorphic to the direct product of 
simply connected, almost simple groups, with $F$-action.  Now it is
easy to see that we are reduced to the case where 
$G \simeq G_1\times \cdots \times G_r$, with $G_i$ a copy of
$G_1$, and $F$ acts on $G$ as a cyclic permutation of all the factors.
Then $G_1$ is $F^r$-stable, and  the set $\CM^F_G$
is in bijective correspondence with the set 
$\CM^{F^r}_{G_1}$, via the correspondence
$(L,C,\CE) \lra (L_1, C_1, \CE_1)$, where
\begin{align*}
L &= L_1 \times F^{-r+1}(L_1)\times\cdots\times F\iv(L_1),  \\
C &= C_1\times F^{-r+1}(C_1)\times\cdots\times F^{-1}(C_1),  \\ 
\CE &= \CE_1\boxtimes 
    {F^{r-1}}^*\CE_1\boxtimes\cdots\boxtimes F^*\CE_1.
\end{align*}
Moreover, $C_1^{F^r} \simeq C^F $ via 
$v_1 \mapsto v = (v_1, F(v_1), \dots, F^{r-1}(v_1))$.
Then $\vf_0: F^*\CE \isom \CE$ is determined by 
$\vf_1: {F^r}^*\CE_1 \isom \CE_1$ as 
\begin{equation*}
(\vf_0)_v = (\vf_1)_{v_1}\otimes(\vf_1)_{F^{r-1}(v_1)}
             \otimes\cdots\otimes(\vf_1)_{F(v_1)}
 \end{equation*}
on $\CE_v = (\CE_1)_{v_1}\otimes(\CE_1)_{F^{r-1}(v_1)}
                 \otimes\cdots\otimes(\CE_1)_{F(v_1)}$.
Similarly, the mixed $F$-structure of $(C',\CE') \in \CN_G$ is described 
by the mixed $F^r$-structure of  $({C'}_1, {\CE'}_1) \in \CN_{G_1}$.
\par
Thus, the determination of the mixed structure of $(C',\CE')$ is 
reduced to the case where $G$ is an $F$-stable, simply connected, 
almost simple group. 
\para{1.6.}
Assume that $G$ is almost simple and simply connected.  Let 
$\Fg = \Lie G$ be
the Lie algebra of $G$.  We further assume that $p$ is good for $G$
unless $G$ is of type $A$, and that $p > n$ if $G = SL_n$.  
Then by [BR], there exists a logarithm map $\log: G \to \Fg$  
satisfying the following properties; $\log$ is an $Ad(G)$-equivariant 
morphism and $\log (1) = 0, d(\log)_1: \Fg \to \Fg$ is the identity
map.  In particular, for any closed subgroup $H$ of $G$, 
$\log(H) \subset \Lie H \subset \Fg$.
Moreover, $\log|_{G\uni}$ turns out to be an isomorphism 
$G\uni \to \Fg\nil$, 
where $\Fg\nil$ is the nilpotent variety of $\Fg$.
\par
Let $\CL$ be an irreducible $G$-local system on a nilpotent orbit
$\ZC$ in $\Fg$.  The notion of cuspidal local system on $\ZC$ is
defined in a similar way as in the case of groups, i.e., $\CL$ is said
to be cuspidal or $(\ZC, \CL)$ is a cuspidal pair 
if for any proper parabolic subalgebra $\Fp_1$ of $\Fg$
with nilpotent radical $\Fn_1$ and any $y \in \Fp_1$, we have
$H_c^i((y + \Fn_1) \cap \ZC, \CL) = 0$ for any $i$.
Then it is easily checked (cf. [L4]) that $\log^*$ gives a bijection 
between the set of cuspidal pairs in $G$ and the set of cuspidal pairs
in $\Fg$.
\par
Let $(L, C, \CE) \in \CM_G$, and $(\ZC, \CL)$ the
corresponding cuspidal pair in  $\Fl = \Lie L$, where 
$C = \log\iv(\ZC), \CE = \log^*\CL$.  We put $\Fp = \Lie(P)$
and $\Fn_P = \Lie U_P$.
Let $\ZC' = \log(C')$ be a nilpotent orbit in $\Fg$.
For each $y \in \ZC'$, put 
\begin{align*}
\tag{1.6.1}
\CP_y &= \{ gP \in G/P \mid \Ad(g)\iv y \in \ZC + \Fn_P\}, \\
\wh \CP_y &= \{ g \in G \mid \Ad(g)\iv y \in \ZC + \Fn_P\}. 
\end{align*} 
Then by using  a similar diagram as in (1.4.2), one can define
a local system $\dot\CL$ on $\CP_y$.  It is easy to see that $\log$
gives an isomorphism $\wh \CP_u \isom \wh \CP_y$ with $y = \log(u)$, 
and so induces an isomorphism $\CP_u \isom \CP_y$.  Then we have 
$\log^*\dot\CL  = \dot\CE$.  It follows that we have a canonical 
isomorphism
\begin{equation*}
\tag{1.6.2}
H_c^{a_0+r}(\CP_u, \dot\CE) \simeq H_c^{a_0+r}(\CP_y, \dot\CL).
\end{equation*}
\par
In the case where $G$ has an $\Fq$-structure with Frobenius map $F$, 
$\Fg$ has also an action of $F$, and we may assume that 
$\log$ is $F$-equivariant.  Then the isomorphism (1.6.2) is compatible 
with $\Fq$-structures.  
We denote by the same symbol $\F$ the linear map on 
$H_c^{a_0+r}(\CP_y, \dot\CL)$ obtained as in the case of 
$H_c^{a_0+r}(\CP_u,\dot\CE)$.   
Hence the linear map $q^{(a_0+r)/2}\p$ on 
$\CE_u$ can be described in terms of the Frobenius action $\F$ on 
$H_c^{a_0+r}(\CP_y, \dot\CL)_{\r}$.
 
\section{Graded Hecke algebras}
\para{2.1.}
The graded Hecke algebra $\BH$ was introduced by Lusztig [L7], which 
is a degenerate version of affine Hecke algebras.  In this section,
following [L7] we review the definition of $\BH$ and their
representations on equivariant $K$-homology groups.
In [L7], $\BH$ are constructed as algebras over $\BC$, but here we
regard them as algebras over $\Ql$ so that one can relate them to 
$l$-adic cohomology groups.
\par  
Let $\F$ be a root system with a set of simple roots 
$\Pi = \{ \a_1, \dots, \a_m\}$ and $W$ the Weyl group of $\F$ with
corresponding simple reflections $\{ s_1, \dots, s_m\}$.   
We assume that the root lattice $\BZ \F$ is embedded in 
a vector space $\Fh^*$ over $\Ql$.  The action of $W$ on $\BZ \F$
makes $\Fh^*$ into a $W$-module.  (Hence $\Fh^*$ has a direct sum
decomposition, one summand being $W$-invariant, the other having 
$\Pi$ as a basis).   
Let $\BS$ be the symmetric algebra of $\Fh^*\oplus \Ql$. We denote
$\Br = (0,1) \in \Fh^*\oplus \Ql$, so that 
$\BS = S(\Fh^*)\otimes \Ql[\Br]$.
$W$ acts naturally on $\BS$ so that $\Br$ is left invariant by $W$.
We denote by $\xi \mapsto {}^w\xi$ the action of $W$ on 
$\CS$.  
Let $c_1, \dots, c_m$ be integers $\ge 2$ such that $c_i = c_j$
whenever $s_i$ and $s_j$ are conjugate in $W$. 
Let $e$ be the neutral element of $W$. 
Lusztig showed in [L7, Theorem 6.3] that there is a unique structure of
associative $\Ql$-algebra on the $\Ql$-vector space 
$\BH = \BS\otimes \Ql[W]$ with unit $1\otimes e$  such that
\par\medskip
(i) \ $\xi \mapsto \xi\otimes e$ is an algebra 
homomorphism $\BS \to \BH$,
\par
(ii) \ $w \mapsto 1\otimes w$ is an algebra homomorphism 
$\Ql[W] \to \BH$, 
\par
(iii) \ ($\xi\otimes e)\cdot (1\otimes w) = \xi\otimes w, \quad
            (\xi \in \BS, w \in W)$,  
\par
(iv) \ $\displaystyle (1\otimes s_i)(\xi\otimes e) - 
     ({}^{s_i}\xi\otimes e)(1\otimes s_i) = 
           c_i\Br\frac{\xi - {}^{s_i}\xi}{\a_i}\otimes e, 
              \quad (\xi \in \BS, 1 \le i \le m)$.
\par\medskip
$\BH$ is called a graded Hecke algebra attached to $W$ with 
parameters $c_i$.
It follows from (iv) we see that $\Br$ is in the center of $\BH$.
\para{2.2.} 
The discussion in [L7] is concerned with algebraic groups over $\BC$.
Hence the equivariant $K$-homology is defined for the varieties over
$\BC$.  Since we treat algebraic groups over finite fields, we need
to construct the equivariant $K$-homology based on the $l$-adic
cohomology groups.  Fortunately, the basic properties established
in section 1 in [L7] work well also for our situation, by a suitable
modification. We give some comments below.
\par  
Let $G$ be an affine algebraic group over $k$, and let $X$ be 
a $k$-variety on which $G$ acts algebraically.  As in [L7], 
for each integer $m \ge 1$, there exists a smooth irreducible 
variety $\vG$ with free $G$-action such that $\vG \to G\bs\vG$ has a locally
trivial principal $G$-fibration, and that  
$H^i(\vG,\Ql) = 0$ for
$i = 1, \dots, m$. 
(As in [L7, 1.1], we embed $G$ as a closed subgroup of $GL_r$, and
consider the embedding 
\begin{equation*}
\tag{2.2.1}
G \subset GL_r \times \{ e\} \subset GL_r\times GL_{r'}
                  \subset GL_{r+r'}.
\end{equation*}
Then $\vG = (\{ e\}\times GL_{r'}) \bs GL_{r+r'}$ for large $r'$
$(2r' \ge m+2)$, with the left action of $G$ on $\vG$, 
satisfies the required condition.) 
For a $G$-variety $X$,  we consider 
${}_{\vG}X = G\bs (\vG \times X)$ 
(the quotient by the diagonal action of $G$).  Then for an 
$G$-equivariant local system $\CL$ 
on $X$, there exists a unique local system ${}_{\vG}\CL$ on
${}_{\vG}X$ such that $\pi^*({}_{\vG}\CL) = p^*\CL$, where
$\pi: \vG \times X \to G\bs (\vG \times X)$
is a natural map, 
and $p: \vG \times X \to X$ is a projection. 
  Then as in [L7], we define
\begin{equation*}
H_G^j(X,\CL) = H^j({}_{\vG}X, {}_{\vG}\CL), \quad
H^G_j(X, \CL) = H_c^{2d-j}({}_{\vG}X, {}_{\vG}\CL^*)^*, 
\end{equation*}      
where $d = \dim ({}_{\vG}X)$, and the upper-script * denotes the 
dual local system or the dual vector space.
(We understand that $H^j_G(X,\CL) = H^j(X, \CL)$ and
$H^G_j(X, \CL) = H_c^{2\dim X - j}(X, \CL^*)^*$ in the case where
$G = \{ e\}$.)  We write them as $H_G^j(X), H_j^G(X)$ if $\CL$ is a 
constant sheaf $\Ql$.  Also we write $H_c^i(X)$, $H^i(X)$ 
instead of $H^i_c(X,\Ql)$, $H^i(X,\Ql)$.
\par
By cup-product, $H_G^*(X) = \bigoplus_j H_G^j(X)$ becomes a 
graded $\Ql$-algebra with 1, and
\begin{equation*}
H_G^*(X,\CL) = \bigoplus_j H^j_G(X,\CL), \quad
H^G_*(X,\CL) = \bigoplus_j H_j^G(X,\CL)
\end{equation*}
become graded $H^*_G(X)$-modules.
\par
We write $H_G^*, H^G_*$ instead of 
$H_G^*(\text{point}), H^G_*(\text{point})$.
Then the map $X \to \text{point}$ defines a $\Ql$-algebra
homomorphism $\ve: H_G^* \to H_G^*(X)$ preserving the
grading.  Via the map $\ve$, $H^*_G(X,\CL), H^G_*(X,\CL)$
can be regarded also as $H_G^*$-modules.
\para{2.3.}
Let $T$ be a torus and $X(T)$ be its character group.
The arguments in 1.10 in [L7] do not hold in that form.  We modify
them as follows.  In the case where $T \simeq \BG_m$ is the 
one dimensional torus,
it can be verified directly by the definition that 
$H^*_T \simeq \Ql[x]$, a polynomial ring with one variable, with 
$x \in H^2_T$. Since 
$H^*_{G\times G'} \simeq H^*_G\otimes H^*_{G'}$, 
we see that $H^*_T \simeq S(V^*)$, the symmetric algebra of
a $\Ql$-vector space $V^* = \Ql\otimes_{\BZ}X(T)$.    
In particular, we have 
\begin{equation*}
H_T^{2j} \simeq S^j(V^*), \quad H_T^{2j+1} = 0, 
\end{equation*}
and we may identify $H^2_T$ with $V^*$.
($S^j(V^*)$ denotes the degree $j$-part of $S(V^*)$).
\par
For $\x \in X(T)$, let $k_{\x}$ be the $T$-module $k$ 
with the $T$-action by $(t, z) \mapsto \x(t)z$.
Let $i: \{ 0\} \hra k, \pi: k \to \{ 0\}$ be the obvious maps.
Then $\pi^*$ is an isomorphism, and the composition
\begin{equation*}
\begin{CD}
H^T_*(\{0\}) @>i_!>> H^T_*(k_{\x}) 
                   @>(\pi^*)\iv >>  H^T_*(\{ 0\})
\end{CD}
\end{equation*}
is $H^*_T$-linear of degree 2.  Since 
$H^T_*(\{0\}) \simeq H_T^*$ as $H^*_T$-modules,  
$(\pi^*)\iv\circ i_!$ is given by multiplication by an element
$c(\x) \in H_T^2$ (cf. [L7, 1.10]).  The map 
$c: X(T) \to H_T^2 = V^*, \x \mapsto c(\x)$ gives an injective 
 group homomorphism.
\par
 Assume that $G$ is an algebraic group such that $G^0$ is a torus 
$T$. Then $W = G/G^0$ acts naturally on $H^*_T$, preserving the grading 
(see [L7, 1.9]). 
$W$ acts also on $X(T)$, and we have
\par\medskip\noindent
(2.3.1) \ The map $c: X(T) \to H_T^2 = V^*$ is $W$-equivariant.
\par\medskip
In fact, take $\vG$ on which $G$ acts freely.  Then, for
a representative $\dw \in G$ of $w \in W$,
 the map
$\vG \times k_{\x} \to \vG \times k_{w(\x)}, (g,x) \mapsto (\dw g, x)$ 
induces a map 
$f_w: T\backslash (\vG \times k_{\x}) \to 
         T\backslash(\vG \times k_{w(\x)})$,
which makes the following diagram commutative.
\begin{equation*}
\begin{CD}
H^T_*(\{0\}) @>i_!>> H^T_*(k_{\x}) 
                   @>(\pi^*)\iv >>  H^T_*(\{ 0\}) \\
@Vw VV                 @VV(f_w^*)\iv V            @VV w V \\
H^T_*(\{0\}) @>i_!>> H^T_*(k_{w(\x)}) 
                   @>(\pi^*)\iv >>  H^T_*(\{ 0\}).
\end{CD}
\end{equation*}
(2.3.1) follows from this.
\par
It follows from (2.3.1) that we have
\begin{equation*}
\tag{2.3.2}
H^*_T \simeq S(\Ql\otimes_{\BZ}X(T)) 
\end{equation*}
as graded $W$-modules. 
\par
We don't know whether the counter part of 1.11 in [L7] holds 
in our setting. However, the following related fact holds.
\begin{lem}
Assume that $G$ is a connected algebraic group.  Let $G_r$ be 
a maximal reductive subgroup of $G$, and $T$ a maximal
torus of $G_r$.  Let $W =  N_{G_r}(T)/T$ be the Weyl group of $G_r$.
Then $W$ acts naturally on $H^*_T$, and the natural map 
$H^*_{G} \to H_T^*$ 
(cf. [L7, 1.4 (g)]) induced from the inclusion $T \hra G$ gives 
an isomorphism 
\begin{equation*}
H^*_G \isom (H^*_T)^W.
\end{equation*} 
\end{lem}
\begin{proof}
By [L7, 1.4 (h)], we know that $H^*_G \isom H^*_{G_r}$.
Hence it is enough to show the lemma in the case where 
$G$ is reductive.  Assume that $G = G_r$.
Let $m$ be a large integer and let $\vG$ be an irreducible, 
smooth variety with a free $G$-action such that $H^i(\vG) = 0$
for $ 1\le i \le m$.  We consider the map
$f:T\backslash\vG \to G\backslash\vG$, which is a locally trivial
fibration with fibre isomorphic to $T\backslash G$. 
We have a spectral sequence
\begin{equation*}
\tag{2.4.1}
H^p(G\bs\vG, R^qf_*\Ql) \Rightarrow 
       H^{p+q}(T\bs \vG). 
\end{equation*}
The map $f$ is $W$-equivariant with respect to the trivial action of 
$W$ on $G\bs\vG$, and the left action of $W$ on $T\bs\vG$, and so
$R^qf_*\Ql$ has a structure of $W$-sheaf, which induces an action 
of $W$ on $H^p(G\bs\vG, R^qf_*\Ql)$.  $W$ acts naturally on 
$H^{p+q}(T\bs\vG)$, and by taking the $W$-invariant parts in 
(2.4.1),  we have a spectral sequence 
\begin{equation*}
\tag{2.4.2}
H^p(G\bs\vG, R^qf_*\Ql)^W \Rightarrow H^{p+q}(T\bs\vG)^W. 
\end{equation*}
Since $f$ is a locally trivial fibration, $R^qf_*\Ql$ is a
local system with fibre $H^q(T\bs\vG)$.  
We may assume that $\vG = (\{e\}\times GL_{r'})\backslash GL_{r+r'}$
as in 2.2.  Then $f$ is $GL_{r+r'}$-equivariant, and so 
$R^qf_*\Ql$ is a $GL_{r+r'}$-local system on the space
$G\bs\vG$ (with respect to the right action of $GL_{r+r'}$).  
Now $GL_{r+r'}$ acts transitively on $G\backslash \vG$ with 
a stabilizer of a point isomorphic to $G \times GL_{r'}$.
Since $G$ is connected, we see that $R^qf_*\Ql$ is
a constant sheaf $H^q(T\bs\vG)$.  It follows that
\begin{equation*}
H^p(G\bs\vG, R^qf_*\Ql) \simeq 
      H^p(G\bs\vG)\otimes H^q(T\bs G) 
\end{equation*}
and we have
\begin{equation*}
H^p(G\bs\vG, R^qf_*\Ql)^W \simeq 
           H^p(G\bs\vG)\otimes H^q(T\bs G)^W
\end{equation*}
since $W$ acts trivially on $H^p(G\bs\vG)$.
It is known that $H^*(T\bs G)$ is a graded regular $W$-module, and
\begin{equation*}
H^q(T\bs G)^W = \begin{cases}
                       \Ql  &\quad\text{ if } q = 0, \\
                        0   &\quad\text{ otherwise }.
                  \end{cases}
\end{equation*}
Hence the spectral sequence (2.4.2) collapses, and we have
\begin{equation*}
\tag{2.4.3}
H^p(G\bs\vG) \simeq H^p(T\bs\vG)^W.
\end{equation*}
This isomorphism is induced from the natural map 
$H^p(G\bs\vG) \to H^p(T\bs\vG)$.
Since $H^p_G = H^p(G\bs\vG)$, and 
$H^p_T = H^p(T\bs\vG)$ by definition, 
the lemma follows from (2.4.3).
\end{proof}
For later discussion, we note the following.
\begin{cor}
Assume that $G$ is connected reductive, and let $T, W$ be as before.  
Let $L$ be a Levi subgroup of a parabolic
subgroup of $G$ containing $T$.  Assume further that 
$L$ contains a cuspidal pair as in 1.1.  
Put $\CW = N_G(Z_L^0)/L = N_G(L)/L$.  Then the image of
the natural map $H^*_G \to H^*_{Z_L^0}$ coincides with 
$(H^*_{Z_L^0})^{\CW}$.
\end{cor}
\begin{proof}
The inclusions $Z_L^0 \hra T \hra G $ induces the maps
$H^*_G \to H^*_T \to H^*_{Z_L^0}$.  
Put $V^* = \Ql\otimes_{\BZ} X(T)$, $V^*_1 = \Ql\otimes_{\BZ} X(Z_L^0)$.
Then by (2.3.2), the map $H^*_T \to H^*_{Z_L^0}$ is nothing but
the natural map $\vf: S(V^*) \to S(V^*_1)$ obtained from the restriction map
$X(T) \to X(Z_L^0)$.  Now $W, \CW$ acts naturally on $S(V^*), S(V^*_1)$,
respectively.  Since $\CW \simeq N_W(W_L)/W_L$, $\vf$ induces a map
$\wt\vf : S(V^*)^W \to S(V^*_1)^{\CW}$.  
By [L7, Proposition 2.6], $Z_L^0$ coincides with a maximal torus of a
certain connected reductive subgroup $H$ of $G$, and $\CW$ is regarded 
as the Weyl group of $H$.  Thus in view of Lemma 2.4, 
it is enough to show that $\wt\vf$ is surjective.
This is equivalent to the fact that $V_1/\CW \to V/W$ is a closed
embedding, where $V$ is the dual space of $V^*$ which is identified
with the Lie algebra of the torus $T_{\Ql}$ over $\Ql$, and similarly 
for $V_1$.  But by using the classification of the triple 
$(L,C,\CE) \in \CM_G$, it is checked that $V_1/\CW \to V/W$ is a
closed embedding.  Thus the corollary follows.
\end{proof}
\para{2.6.}
1.12 (a), (b) in [L7] were deduced by using 1.11 there.  
Here we show the corresponding facts by using 2.3 as follows. 
\par\medskip\noindent
(2.6.1) \ Let $G$ be an algebraic group such that $G^0$ is a central 
torus in $G$.  Then we have
\begin{equation*}
H^*_G \simeq H^*_{G^0}.
\end{equation*} 
\par\medskip
In fact, by [L7, 1.9 (a)], we have
\begin{equation*}
H^*_G \simeq (H^*_{G^0})^{G/G^0}.
\end{equation*} 
But $H^*_{G^0} \simeq S(V^*)$ with $V^* = H^2_{G^0}$, and the action of 
$G/G^0$ on $S(V^*)$ is determined by the action of $G/G^0$ on $X(G^0)$
by (2.3.2).  By our assumption, $G/G^0$ acts trivially on $X(G^0)$, 
and so on $S(V^*)$.  This implies that 
$H^*_G \simeq S(V^*) \simeq H^*_{G^0}$, and (2.6.1) follows. 
\par\medskip\noindent
(2.6.2) \ In the same setting as above, let $E$ be an irreducible
representation of $G/G^0$ over $\Ql$.
Then we have 
\begin{equation*}
H^G_*(\text{point}, E\otimes E^*) \simeq H^{G^0}_*.
\end{equation*} 
\par\medskip
The proof is similar to [L7, 1.12 (b)], by making use of (2.6.1).
\para{2.7.}
We return to the setting in 1.1, and consider a connected reductive 
algebraic group $G$, and its Lie algebra $\Fg$. 
We further assume that $G$ is almost simple, simply connected.
Let $\BG_m$ be the
multiplicative group of $k$.  Then $G$ acts on $\Fg$
by the adjoint action, and $G \times \BG_m$ acts on $\Fg$ by
$(g_1, t): x \mapsto t^{-2}\Ad(g_1)x$.  For $x \in \Fg$, we denote
by $Z_G(x)$ the stabilizer of $x$ in $G$, and by $M_G(x)$ the
stabilizer of $x$ in $G \times \BG_m$.  Hence
\begin{equation*}
M_G(x) = \{ (g_1, t) \in G \times \BG_m \mid Ad(g_1)x = t^2x\}.
\end{equation*}
\par
We assume that $p$ is large enough so that Jacobson-Morozov's theorem
and Dynkin-Kostant theory hold for $\Fg$, (e.g., $p > 3(h-1)$, where
$h$ is the Coxeter number of $W$, [C, 5.5]).  Then, for each nilpotent
element $y \in \Fg$, there exists a Lie algebra homomorphism 
$\f: \Fs\Fl_2 \to \Fg$, and elements $y^-, h \in \Fg$ such that
\begin{equation*}
y = \f\begin{pmatrix}
         0 & 1 \\
         0 & 0
       \end{pmatrix},
\quad
y^- = \f\begin{pmatrix}
          0 & 0 \\
          1 & 0
        \end{pmatrix},
\quad
h = \f\begin{pmatrix}
          1 & 0  \\
          0 & -1 
       \end{pmatrix}.
\end{equation*}   
Thus we have $[h, y] = 2y, [h, y^-] = -2y^-, [y, y^-] = h$.
Moreover, we have a decomposition $\Fg = \bigoplus_i \Fg_i$, 
where $\Fg_i$ is the $i$-eigenspace of $\ad h: \Fg \to \Fg$. 
In particular, note that $y \in \Fg_2, y^- \in \Fg_{-2}$.
One can define an algebra homomorphism $\r': \BG_m \to \Aut \Fg$
by $\r'(t)z = t^iz$ for $z \in \Fg_i$.  Since the identity 
component of $\Aut \Fg$ coincides with $\ad G = G/Z_G$,  
$\r'(\BG_m)$ is a one-dimensional torus in $\ad G$. By taking 
the identity component of $\pi\iv(\r'(\BG_m))$ for $\pi: G \to \ad G$,
one obtains a one parameter subgroup $\r : \BG_m \to G$ such that
$\r' = \pi\circ\r$.
\par
We put 
\begin{align*}
Z_G(\f) &= Z_G(y) \cap Z_G(y^-), \\
M_G(\f) &= \{ (g_1, t) \in G \times \BG_m \mid
  \Ad(g_1)y = t^2y, \Ad(g_1)y^- = t^{-2}y^- \}.
\end{align*}  
It is known that $Z_G(\f)$ is a maximal reductive subgroup of 
$Z_G(y)$.  
It is easy to check that $(g_1, t) \mapsto (g_1\r(t), t)$
gives isomorphisms of algebraic groups
\begin{equation*}
\tag{2.7.1}
Z_G(y)\times \BG_m \isom M_G(y), \qquad
Z_G(\f)\times \BG_m \isom M_G(\f).
\end{equation*}
Hence $M_G(\f)$ is also a maximal reductive subgroup of $M_G(y)$.
It also follows from (2.7.1) that the embedding  
$Z_G(y) \hra M_G(y)$ by $g_1 \mapsto (g_1, 1)$ induces 
an isomorphism
\begin{equation*}
Z_G(y)/Z_G^0(y) \isom M_G(y)/M_G^0(y).
\end{equation*} 
This implies that the $G$-orbit of $x \in \Fg$ is also a $G\times \BG_m$-orbit,
and a $G$-local system on a nilpotent $G$-orbit in $\Fg$ is
automatically a $G \times \BG_m$-local system.
In later discussions, we use the notation $M(y), M^0(y)$, etc. instead
of $M_G(y), M_G^0(y)$, etc. by omitting the subscript $G$ if there is no
fear of confusion. 
\para{2.8.}
Under the setting in 1.1, let  $\Fp, \Fl, \Fn_P$ be the 
Lie algebras of $P, L, U_P$ so that 
$\Fp = \Fl \oplus \Fn_P$.
Let $\Fz$ the Lie algebra of $Z_L^0$. 
We assume that $(L, C, \CE) \in \CM_G$, and let $(\ZC, \CL)$
be the corresponding cuspidal pair on $\Fl$ (cf. 1.6.)
Let
\begin{equation*}
\tag{2.8.1}
\dot\Fg = \{ (x, gP) \in \Fg \times G/P \mid 
            \Ad(g\iv)x \in \ZC + \Fz + \Fn_P \},
\end{equation*}
and $\pi: \dot\Fg \to \Fg$ be the first projection.   
$G \times \BG_m$ acts on $\dot\Fg$ by 
$(g_1, t): (x, gP) \mapsto (t^{-2}\Ad(g_1)x, g_1gP)$,
and $\pi$ is $G \times \BG_m$-equivariant.
We consider the diagram
\begin{equation*}
\begin{CD}
\ZC @<\a << 
 \wh{\dot\Fg} = \{ (x, g) \in \Fg \times G \mid 
         \Ad (g\iv)x \in \ZC + \Fz + \Fn_P  \} @>\b>> \dot\Fg, 
\end{CD}
\end{equation*}
where $\a(x, g) = \pr_{\ZC}(\Ad(g\iv)x), \b(x,g) = (x, gP)$.
Here $\a, \b$ are $G \times \BG_m$-equivariant with respect to 
the action of  $G \times \BG_m$ on $\ZC$ given by 
$(g_1, t): x \mapsto t^{-2}x$,  
and the action of it on the middle term given by 
$(g_1, t): (x,g) \mapsto (t^{-2}\Ad(g_1)x, g_1g)$.
Since $\CL$ is an $L$-local system, there exists a unique
local system $\dot\CL$ on $\dot\Fg$ such that 
$\a^*\CL = \b^*\dot\CL$.  By 2.7, $\CL$ is  
$L \times \BG_m$-equivariant, and so is $G \times \BG_m$-equivariant
with respect to the above action.  Hence $\dot\CL$ turns out to be 
$G \times \BG_m$-equivariant. 
\par
Let $\dot\CL^*$ be the dual local system of $\dot\CL$, and 
consider $K = \pi_!(\dot\CL^*)$.  Then it is shown in [L7, 3.4]
that $K[\d]$ is a $G \times \BG_m$-equivariant perverse sheaf on $\Fg$
with a canonical $\CW$ action, where 
$\d = \dim (\Fg/\Fl) + \dim (\ZC + \Fz)$.
\par
Let $X$ be an algebraic variety with a given morphism $m: X \to \Fg$.
We consider the fibre product $\dot X = X\times_{\Fg}\dot\Fg$ with
the cartesian diagram
\begin{equation*}
\begin{CD}
\dot X @>\dot m >> \dot \Fg \\
@V \pi' VV            @V\pi VV \\
X @> m >>           \Fg
\end{CD}
\tag{2.8.2}
\end{equation*}
Then $m^*K$ is a complex with $\CW$-action, and it induces
a natural $\CW$-action on the cohomologies 
\begin{equation*}
\tag{2.8.3}
\HH^j_c(X, m^*K) \simeq \HH^j_c(X, \pi'_!\dot m^*\dot\CL^*)
                   \simeq H_c^j(\dot X, \dot m^*\dot\CL^*). 
\end{equation*}
\par
We further assume that $X$ is a $G'$-variety, where $G'$ is a
connected closed subgroup of $G \times \BG_m$, and that
$m$ is compatible with $G'$-actions.
If we choose a smooth irreducible variety $\vG$ with a free 
$G'$-action as in 2.2, 
the cartesian diagram (2.8.2) is lifted to the cartesian diagram
\begin{equation*}
\begin{CD}
{}_{\vG}\dot X @>{}_{\vG}\dot m>> {}_{\vG}\dot\Fg \\
@V {}_{\vG}\pi' VV            @V{}_{\vG}\pi VV \\
{}_{\vG}X @> {}_{\vG}m >>           {}_{\vG}\Fg
\end{CD}
\end{equation*}
As in 2.2, we have a local system ${}_{\vG}\dot\CL^*$ on 
${}_{\vG}\dot\Fg$, and a perverse sheaf (up to shift) ${}_{\vG}K$
on ${}_{\vG}\Fg$ which inherits a $\CW$-action from $K$.
Since ${}_{\vG}K = ({}_{\vG}\pi)_!({}_{\vG}\dot\CL^*)$, as in 
(2.8.3) we have natural $\CW$-actions on cohomologies
\begin{equation*}
\HH^j_c({}_{\vG}X, ({}_{\vG}m)^*({}_{\vG}K)) 
\simeq 
\HH^j_c({}_{\vG}X, ({}_{\vG}\pi')_!({}_{\vG}
        \dot m)^*{}_{\vG}\dot\CL^*)
\simeq
H^j_c({}_{\vG}\dot X, ({}_{\vG}\dot m)^*{}_{\vG}\dot\CL^*).
\end{equation*}
Hence we have an action of $\CW$ on the equivariant homology 
\begin{equation*}
H_j^{G'}(\dot X, \dot\CL) = 
   H_c^{2d-j}({}_{\vG}\dot X, {}_{\vG}\dot\CL^*)^*, 
\end{equation*}
where $d = \dim ({}_{\vG}\dot X)$.
(Here we write  $\dot m^*\dot\CL^*$, 
$({}_{\vG}\dot m)^*{}_{\vG}\dot\CL^*$, etc. as $\dot\CL^*$, 
${}_{\vG}\dot\CL^*$, etc. by abbreviation.) 
\para{2.9.}
We fix an element $x_0 \in \ZC$ and a Lie algebra homomorphism
$\f_0: \Fs\Fl_2 \to \Fl$ such that 
$\f_0\binom {0 1}{0 0} = x_0$.
As in [L7, 2.3 (b)], we have
\begin{equation*}
\tag{2.9.1}
Z_L^0(\f_0) = Z_L^0.
\end{equation*}
It follows that $Z^0_L(\f_0)$ is central in $Z_L(\f)$.  Hence
by (2.7.1), we see that 
\par\medskip\noindent
(2.9.2)
\ $M_L^0(\f_0) \simeq Z_L^0 \times \BG_m$, and 
$M_L^0(\f_0)$ is contained in the center of $M_L(\f)$.
\par\medskip
Put $\Fh^* = \Ql\otimes_{\ZZ}X(Z_L^0)$.
The $\Fh^*$ is a $\Ql$-space of $\dim_{\Ql} \Fh^* = \dim_k \Fz$, 
on which $\CW$ acts naturally. We define a symmetric algebra 
$\BS$ over $\Ql$ by 
\begin{equation*}
\BS = S(\Fh^*\oplus\Ql) = S(\Fh^*)\otimes \Ql[\Br],
\end{equation*}
where $\Ql[\Br]$ is the polynomial ring with an indeterminate 
$\Br$ corresponding to $(0,1) \in \Fh^*\oplus \Ql$.
We now consider the equivariant cohomology 
$H^*_{G\times \BG_m}(\dot\Fg)$.
As in [L7, Proposition 4.2],  we have an isomorphism 
\begin{equation*}
\tag{2.9.3}
H^*_{G\times \BG_m}(\dot\Fg) \simeq \BS
\end{equation*}
as graded algebras.  In particular, $H^j_{G\times\BG_m}(\dot\Fg) = 0$
for odd $j$.  For the proof, the argument in [L7] implies that
\begin{equation*}
H^*_{G\times\BG_m}(\dot\Fg) \simeq H^*_{M_L(\f_0)}.
\end{equation*}
Then by using (2.6.1) and (2.9.2), combined with (2.3.2), we have 
\begin{equation*}
H^*_{M_L(\f_0)} \simeq H^*_{M_L^0(\f_0)} \simeq 
                    H^*_{Z_L^0\times\BG_m} = \BS. 
\end{equation*}
Hence (2.9.3) follows.
\par
Let $\wt X$ be a $G'$-variety ($G'$ is a connected closed subgroup of 
$G\times\BG_m$), with a given $G'$-equivariant morphism 
$\wt m: \wt X \to \dot\Fg$.
$\wt m^*\dot\CL$ is a $G'$-local system on $\wt X$, which
we denote by $\dot\CL$ by abbreviation.
Now $\wt m^*$ induces an algebra homomorphism 
$H^*_{G'}(\dot\Fg) \to H^*_{G'}(\wt X)$.
By combining the natural homomorphism 
$H^*_{G\times\BG_m}(\dot\Fg) \to H^*_{G'}(\dot\Fg)$ (cf. [L7, 1.4 (g)]),
we have a homomorphism 
$H^*_{G\times\BG_m}(\dot\Fg) \to H^*_{G'}(\wt X)$.
Since $H^{G'}_*(\wt X, \dot\CL)$ is a 
$H^*_{G'}(\wt X)$-module by 2.2, 
$H^{G'}_*(X, \dot\CL)$ has a structure of a left 
$H^*_{G\times \BG_m}(\dot\Fg)$-module.  Thus by (2.9.3), 
$H^{G'}_*(\wt X, \dot\CL)$ turns out to be an $\BS$-module. 
\para{2.10.}
Let $\pi: \dot\Fg \to \Fg$ be as in 2.8.  Then for each $y \in \Fg\nil$, 
$\pi\iv(y)$ coincides with $\CP_y$ in (1.6.1).  The variety
$X = \{ y\}$ is invariant under the action of 
$M^0(y) \subset G \times \BG_m$.
Let $G'$ be a connected closed subgroup of $M^0(y)$.  
By applying 2.8 to the inclusion $m: X \hra \Fg$ together with 
$\dot X = \CP_y$, we see that 
$H^{G'}_*(\CP_y, \dot\CL)$ has a natural $\CW$-action. 
By applying 2.9 for $\wt X = \dot X$, 
$H^{G'}_*(\CP_y,  \dot\CL)$ has a natural $\BS$-action.
It also has a structure of $H^*_{G'}$-module by 2.2.
\par
We consider the graded Hecke algebra $\BH = \BS\otimes\Ql[\CW]$
as defined in 2.1, where $\BS$ is as in 2.9, with a natural action 
of the Coxeter group $\CW$. 
Lusztig proved the following theorem.
\begin{thm}[Lusztig {[L7, Theorem 8.13]}] 
There is a unique $\BH$-module structure on 
$H^{M^0(y)}_*(\CP_y, \dot\CL)$ such that the actions of 
$\BS$ and $\CW$ are given as in 2.10. 
(The integers $c_i$ are determined according to the cuspidal pair
$(\ZC, \CL)$.  See [L7, 2.13] for explicit values for $c_i$).
Moreover, the $\BH$-module
structure commutes with the $H^*_{M^0(y)}$-module structure on
$H^{M^0(y)}_*(\CP_y, \dot\CL)$.
\end{thm}
\par\medskip\noindent
{\bf Remark 2.12.}
\addtocounter{thm}{1}
The arguments used in [L7] to prove the theorem are valid also for
our setting in almost all cases, by taking 2.3 $\sim$ 2.7 into account.
We give further comments on the discrepancies of the arguments.
\par
(a) In [L7, 4.3], the property of the image 
$H^*_{G\times\BG_m} \to H^*_{M_L^0(\f_0)}$ is used. For this
we appeal to Corollary 2.5.
\par
(b) In the proof of Proposition 7.2 in [L7], a property of simply
connected space is used, which is not valid in the positive characteristic
case.  As in 7.1, we consider a connected algebraic group $M$, and an
$M$-variety $X$, $M$-equivariant local system $\CE$ on $X$. 
Let $\vG$ be an irreducible, smooth variety with a free $M$-action
as before.  Let $f: M\bs(\vG\times X) \to M\bs \vG$ be the locally
trivial fibration.  We consider the Leray-Serre spectral sequence
\begin{equation*}
H_c^p(M\bs\vG, R^qf_!({}_{\vG}\CE^*)) \Rightarrow 
      H^{p+q}_c(M\bs(\vG\times X), {}_{\vG}\CE^*).
\end{equation*}
We show that 
\begin{equation*}
\tag{2.12.1}
E_2^{p,q} = H_c^p(M\bs\vG, R^qf_!({}_{\vG}\CE^*)) 
              = H_c^p(M\bs\vG)\otimes H_c^q(X, \CE^*).
\end{equation*}
(In [L7], this is obtained as a consequence of the fact that 
$M\bs\vG$ can be chosen to be simply connected).
We consider the cartesian diagram
\begin{equation*}
\begin{CD}
\vG\times X @>\pi >>  M\bs(\vG \times X) \\
@V\wt f VV               @VVf V          \\
\vG   @>\wt\pi >>        M\bs\vG.
\end{CD}
\end{equation*}
Now ${}_{\vG}\CE^*$ on $M\bs(\vG\times X)$ satisfies the property that
$ \Ql\boxtimes \CE^* = \pi^*({}_{\vG}\CE^*)$.
By the base change theorem, we have 
$\wt\pi^*R^qf_!({}_{\vG}\CE^*) \simeq R^q\wt f_!\pi^*({}_{\vG}\CE^*)$.
It is easy to see that $R^q\wt f_!(\Ql\boxtimes \CE^*)$ is an 
$M$-equivariant  constant sheaf, and $R^qf_!({}_{\vG}\CE^*)$ is 
obtained from it as the unique quotient.  Thus, $R^qf_!({}_{\vG}\CE^*)$
is also a constant sheaf with the stalk $H^q_c(X,\CE^*)$. 
This implies (2.12.1).
\par
Once this is established, the other parts in the proof of 
Proposition 7.2 work without change.
\para{2.13.} 
We return to the setting in 2.10. 
Let $T(y)$ be a maximal torus of $M^0(y)$ and $W(y)$ the Weyl group
of a maximal reductive subgroup of $M^0(y)$ with respect to $T(y)$. 
Then by (2.3.2) and Lemma 2.4, 
$H^*_{M^0(y)}$ can be identified with $S(V^*)^{W(y)}$, 
where $V^* = \Ql\otimes_{\BZ}X(T(y))$. 
Hence $H^*_{M^0(y)}$ may be regarded as the coordinate ring of 
an affine algebraic variety (over $\Ql$) $V_1 = V/W(y)$, 
where $V$ is the dual space of $V^*$.  Then for each $v \in V_1$, one obtains an
algebra homomorphism $H^*_{M^0(y)} \to \Ql$, $f \mapsto f(v)$.  We denote 
the thus obtained $H^*_{M^0(y)}$-module $\Ql$ by $(\Ql)_v$.
It is known by [L7, 8.6] that $H^{M^0(y)}_*(\CP_y, \dot\CL)$ is a
finitely generated projective $H^*_{M^0(y)}$-module.  It follows that 
$H^{M^0(y)}_*(\CP_y, \dot\CL)$ may be regarded as a space of sections
of algebraic vector bundle $E$ over $V_1$, where the fibre of $E$ at
$v \in V_1$ is given by 
\begin{equation*}
\tag{2.13.1}
E_v = (\Ql)_v\otimes_{H^*_{M^0(y)}}H^{M^0(y)}_*(\CP_y, \dot\CL).
\end{equation*} 
Put $\ol{M}(y) = M_(y)/M^0(y)$.  Then the finite group
$\ol{M}(y)$ acts on $H^*_{M^0(y)}$ as a $\Ql$-algebra automorphism, 
and acts on $H^{M^0(y)}_*(\CP_y, \dot\CL)$ compatible with 
the action of $H^*_{M^0(y)}$.  Also this action of $\ol{M}(y)$ on
$H^{M^0(y)}_*(\CP_y,\dot\CL)$ commutes with the action of 
$\BH$. The action of $\ol{M}(y)$ on
$H^*_{M^0(y)}$ induces an action of $\ol{M}(y)$ on $V_1$, and $E$
turns out to be an $\ol{M}(y)$-equivariant vector bundle over $V_1$.
For each $v \in V_1$, we denote by $\ol{M}(y,v)$ the stabilizer of
$v$ in $\ol{M}(y)$.  Then $\ol{M}(y,v)$ acts naturally on $E_v$. 
\par
Let $\ol{M}(y,v)\wg$ be the set of irreducible representations of 
$\ol{M}(y,v)$ up to isomorphisms.  For each 
$\r \in \ol{M}(y,v)\wg$, put 
$E_{v,\r} = (\r^*\otimes E_v)^{\ol{M}(y,v)}$, where $\r^*$ is the 
dual representation of $\r$.  Then $E_{v, \r}$ is an $\BH$-module, 
and $E_v$ is decomposed as 
\begin{equation*}
E_v = \bigoplus_{\r \in \ol{M}(y,v)\wg}\r\otimes E_{v,\r}.
\end{equation*}
\par
The action of $M(y)$ on $\CP_y, \dot\CL, \dot\CL^*$ induces an
action of $\ol M(y)$ on 
$H_c^*(\CP_y,\dot\CL)$, $H_c^*(\CP_y,\dot\CL^*)$, hence on 
$H_*^{\{ e\}}(\CP_y,\dot\CL) = H_c^*(\CP_y,\dot\CL^*)^*$.
It is known by [L7, 8.10] that $E_{v,\r} \ne 0$ if and only if
$\r$ occurs in the restriction of $\ol{M}(y)$-module 
$H^{\{ e\}}_*(\CP_y,\dot\CL)$ to $\ol{M}(y,v)$.
The $\BH$-modules $E_{v,\r}$ are called standard modules.
\par\medskip\noindent
{\bf Remarks 2.14.}
\addtocounter{thm}{1}
(i) \ Standard modules $E_{v,\r}$ are parametrized in [L7] (i.e., 
in the setting that $G$ and $\Fg$ are defined over $\BC$) 
as $E_{h, r_0, \r}$ in terms of the pair
$(h, r_0) \in \Fg \oplus \BC$ such that $[h, y] = 2r_0y$ 
with $h$ semisimple.  This is also possible in our situation, though
we cannot use the Lie algebra $\Fg$ over $k$.
Since $p$ is good, we have corresponding objects $G_{\BC}, \Fg_{\BC}$,
and the parametrization of nilpotent orbits and the structure of 
$\ol{M}(y)$ are the same for $\Fg_{\BC}$ also.
If we consider the maximal torus $T(y)_{\BC}$ 
in $M(y)_{\BC}$ corresponding to $T(y)$ in $M(y)$, the space $V^*$ may be
identified (under a choice of an isomorphism $\Ql\simeq \BC$) with
the dual of the Cartan subalgebra $\Fh(y)_{\BC}$ of a maximal reductive 
subalgebra $\Fm(y)_{\BC,r}$ of $\Fm(y)_{\BC} = \Lie M(y)_{\BC}$ 
with the action of $W(y)$.
Then the action of $\ol{M}(y)$ on $S(V^*)^{W(y)}$ 
coincides with the
action of $\ol{M}(y)$ on 
$S(\Fh(y)_{\BC}^*)^{W(y)}
 \simeq S(\Fm(y)_{\BC,r}^*)^{M^0(y)_{\BC}}$.
Here
\begin{equation*}
\Fm(y)_{\BC} = \Lie M^0(y)_{\BC} = 
  \{ (x, r_0) \in \Fg_{\BC}\oplus \BC \mid [x, y] = 2r_0y\}. 
\end{equation*}
Moreover, the action of $\ol M(y)$ on $S(\Fm(y)^*_{\BC,r})$ is induced from 
the action of 
$M(y)_{\BC}$, $(g_1,t): (x,r_0) \mapsto (t^{-2}\Ad(g_1)x, t^{-2}r_0)$. 
Hence $V_1$ is identified with the set of semisimple $M^0(y)_{\BC}$-orbits on 
$\Fm(y)_{\BC}$.
This implies, in our case, that $E_{v,\r}$ may be expressed as 
$E_{h, r_0, \r}$,
and $\ol{M}(y,v)$ as $\ol{M}(y, h,r_0)$, if $(h,r_0)$ is a semisimple 
orbit in $\Fg_{\BC}\oplus \BC$ corresponding to $v \in V_1$.
\par
(ii) \ Standard modules play a crucial role in the representation
theory of $\BH$.  The structure of $\BH$-module $E_{v,\r}$ was 
studied throughly in [L8], [L9].  However, the result in [L7] 
is enough for our purpose.  
\par\medskip
In view of the above remarks, the following result of Lusztig can 
be applied to our setting.
\begin{thm}[{[L7, Theorem 8.17]}] 
Let $(h, r_0) \in \Fg_{\BC}\oplus\BC$ be a semisimple element 
such that $r_0 \ne 0$.  Then 
\begin{enumerate}
\item
Let 
$Y_{(h,r_0)} = \{ x \in \Fg_{\BC} \mid [h,x] = 2r_0x\}$.
Then $Y_{(h,r_0)}$ consists of
nilpotent elements, and $Z_{G_{\BC}}(h)$ acts (by the adjoint
action ) on $Y_{(h,r_0)}$ with finitely many orbits. 
\item
Let $y$ be an element in the unique open dense orbit in $Y_{(h,r_0)}$. 
Then $(h,r_0) \in \Fm(y)_{\BC}$.  Let $\r \in \ol{M}(y,h,r_0)\wg$ be such
that $E_{h,r_0,\r} \ne \{ 0\}$.  Then $E_{h,r_0,\r}$ is a simple 
$\BH$-module.
\end{enumerate}
\end{thm}
\para{2.16.}
Here we summarize the properties connecting the equivariant homology 
with the ordinary cohomology.  Let $M$ be a
connected algebraic group, $X$ an $M$-variety and $\CE$ an
$M$-equivariant local system on $X$.  We consider $H^M_*(X,\CE)$.
For each $i$, we define $F^i$ as the $H^*_M$-submodule of 
$H^M_*(X,\CE)$ generated by $\bigoplus_{j \le i}H^M_j(X,\CE)$. 
Then $F^i$ gives a filtration 
$F^0 \subseteq F^1 \subseteq \cdots$ and $F^i = 0$ for $i < 0$.
Put $\Pi_i = H^M_i(X,\CE)/H^M_i(X,\CE) \cap F^{i-1}$.
We have a natural injection  $\Pi_i \to F^i/F^{i-1}$ 
as $\Ql$-spaces.
Since $F^i/F^{i-1}$ is an $H^*_M$-module, this is extended to an  
$H^*_M$-linear map
\begin{equation*}
\tag{2.16.1}
H^*_M\otimes_{\Ql}\Pi_i \to F^i/F^{i-1}.
\end{equation*}
The natural homomorphism $H^M_i(X,\CE) \to H_i^{\{e\}}(X,\CE)$
is zero on $H^M_i(X,\CE) \cap F^{i-1}$, and it factors through
a $\Ql$-linear map 
\begin{equation*}
\tag{2.16.2}
\Pi_i \to H^{\{e\}}_i(X,\CE).
\end{equation*} 
Lusztig showed in  [L7, 7.2] that 
the maps (2.16.1) and (2.16.2) are isomorphisms whenever 
$H^{\odd}_c(X,\CE) = 0$, and in that case we obtain
an isomorphism
\begin{equation*}
\tag{2.16.3}
H^*_M\otimes_{\Ql} H^{\{e\}}_i(X,\CE) \isom F^i/F^{i-1}.
\end{equation*}
\par
We now consider the case where $X = \CP_y$, $\CE = \dot\CL$ 
and $M = M^0(y)$.
It is known that $H^{\odd}_c(\CP_y, \dot\CL) = 0$ by [L3, V, 24.8],
and so the previous argument can be applied.
We consider $E_v$ as in (2.13.1) and $H^*_{M^0(y)}$-module $(\Ql)_v$.  
We define an $\Ql$-space $F^i_v$ by 
$F^i_v = (\Ql)_v\otimes_{H^*_{M^0(y)}}F^i$.  Then 
$F^i_v$ is naturally identified with a quotient of 
$\bigoplus_{j \le i}H^{M^0(y)}_j(\CP_y,\dot\CL)$.
We denote by $f_i: F^{i-1}_v \to F^i_v$ the natural map 
induced from $F^{i-1} \hra F^i$.
It follows from (2.16.3) we have an exact sequence of $\Ql$-spaces
\begin{equation*}
\begin{CD}
F^{i-1}_v @>f_i>> F^i_v @>>> H^{\{e\}}_i(\CP_y, \dot\CL) @>>> 0.
\end{CD}
\tag{2.16.4}
\end{equation*}
In particular, we have
\begin{equation*}
\tag{2.16.5}
F^0_v \simeq H^{\{ e\}}_0(\CP_y, \dot\CL).
\end{equation*}
\para{2.17.}
We consider the $\Fq$-structure on the equivariant homology.
Assume that $G$ and $X$ are defined over $\Fq$ with Frobenius map $F$,
and $G$ acts on $X$ over $\Fq$.  Let $\CE$ be an $G$-equivariant local
system on $X$ such that $F^*\CE \simeq \CE$.  We fix an isomorphism
$\vf: F^*\CE \isom \CE$.  Then  $\vf$ induces natural linear 
isomorphisms on $H^G_*(X, \CE), H^*_G(X, \CE)$, etc.  In fact, 
one can choose a $G$-variety $\vG$ so that $\vG$ is defined over
$\Fq$.  (We may assume that $G$ is an $F$-stable closed subgroup of 
some $GL_r$.  The case where $GL_r$ has a split $\Fq$-structure, 
the construction of $\vG$ in 2.2 works well.
If $GL_r$ is of non-split type, we choose $F = \s_0F_0$,
where $F_0$ is a split Frobenius, and $\s_0$ is an automorphism of
$GL_r$ defined by $\s_0(g) = {}^tg\iv$.  By choosing similar Frobenius
maps for $GL_{r'}$ and $GL_{r+r'}$, the inclusions in (2.2.1)
 are $F$-equivariant.  Hence 
$\vG = \{ e\}\times GL_{r'}\bs GL_{r+r'}$ is defined over $\Fq$.)
Then the maps $\pi: \vG\times X \to {}_{\vG}X$,
$p: \vG \times X \to X$ are defined over $\Fq$.  Hence ${}_{\vG}\CE$
inherits an $\Fq$-structure of $\CE$, which induces a linear map on 
$H_G^j(X,\CE) = H^j({}_{\vG}X, {}_{\vG}\CE)$.  The thus obtained linear
map is independent of the choice of $\vG$.  In fact, if $\vG'$ is
another choice, we have an isomorphism 
$H^j({}_{\vG}X, {}_{\vG}\CE) \isom 
    H^j_{\vG\times\vG'}({}_{\vG\times\vG'}X, {}_{\vG\times\vG'}\CE)$,
etc. as in [L7, 1.1], which are compatible with the induced 
$F$-actions on them.

\par\medskip
\section{$G = SL_n$ with $F$ of split type}
\para{3.1.}
In this section, we assume that $p$ is arbitrary, and consider 
$G = SL_n$ with the standard Frobenius map $F$ on $G$, i.e., for 
$g = (g_{ij}) \in G$, $F(g) = (g^q_{ij})$.  Let $V = k^n$ with the
standard basis $e_1, \dots, e_n$ and we identify $SL_n$ with $SL(V)$. 
\par
Let $\Fg = \Fs\Fl_n$ be the Lie algebra of $G$, and we denote by $F$ the 
corresponding Frobenius map on $\Fg$.  
The unipotent classes in $G$ and nilpotent orbits in $\Fg$ are
parametrized by partitions of $n$, via Jordan normal form.  Let 
$\la = (\la_1,\la_2,\dots, \la_r)$ be a partition of $n$, and let 
$C_{\la}$ (resp. $\ZC_{\la}$ ) be the corresponding unipotent class
in $G$ (resp. nilpotent orbit in $\Fg$).  
Each $\ZC_{\la}$ is $F$-stable, and we construct a specific nilpotent
transformation $y = y_{\la} \in \ZC_{\la}^F$ by defining a basis 
$\{ y^af_j \mid 1 \le j \le r, 0 \le a < \la_j\}$ of $V$ 
obtained from the standard basis as follows;
\begin{equation*}
\tag{3.1.1}
y^af_j = e_i \qquad \text{ with } \quad i = \la_1 + \cdots \la_{j-1} + a.
\end{equation*}
Then $u_{\la} = y_{\la} + 1 \in C_{\la}^F$.
The element $y_{\la} \in \ZC_{\la}^F$ (resp. 
$u_{\la} \in C_{\la}^F$) is called the split element 
corresponding to $\la$.
\para{3.2.}
By [L2], [LS], the generalized Springer correspondence for the case
where $G = SL_n$ is described as follows.  Let $n'$ be the largest divisor of
$n$ which is prime to $p$. Then the center $Z_G$ is a cyclic group of 
order $n'$.  For a divisor $d$ of $n'$, 
consider a Levi subgroup $L$ of $P$ of the type
$A_{d-1}\times \cdots\times A_{d-1}$ ($n/d$-factors).
Let $C$ be the regular unipotent class in $L$.  Then for $v \in C$, 
$A_L(v) = Z_L/Z_L^0 \simeq \BZ/d\BZ$.  Let $\CE$ be an $L$-equivariant local
system on $C$ corresponding to a character $\r_0$ of $A_L(v)$ of order
$d$.  Then $(C, \CE)$ is a cuspidal pair on $L$, and any cuspidal pair
on a Levi subgroup of a parabolic subgroup of $G$ 
is obtained in this way.  Hence for a Levi subgroup $L$ determined 
by $d$, there exist exactly $\vf(d)$ cuspidal pairs in $L$, where 
$\vf$ is the Euler function.
\par
Let $K$ be as in (1.2.2) with respect to the cuspidal pair $(C, \CE)$
on $L$.  Let $C'$ be a unipotent class in $G$
corresponding to a partition $\la = (\la_1, \la_2, \dots, \la_r)$.  Then 
for $u \in C'$, $A_G(u)$ is a cyclic group of order $n'_{\la}$, where
$n'_{\la}$ is the greatest common divisor of $n', \la_1,\la_2, \dots, \la_r$. 
Let $\CE'$ be the local system on $C'$ corresponding to $\r \in A_G(u)\wg$.
The condition for $C'$ such that $\IC(\ol{C'},\CE')$
is a component of $K$ (up to shift) is that each $\la_i$ is divisible
by $d$. In this case $n_{\la}'$ is divisible by $d$, 
and we have a surjective homomorphism $A_G(u) \to A_L(v)$ which
factors through the natural maps $Z_G \to A_G(u)$ and $Z_G \to A_L(v)$.  
Let $\r \in A_G(u)\wg$ be 
the character obtained as the pull back of $\r_0 \in A_L(v)\wg$.
Then $\IC(\ol{C'},\CE')$ is the unique component in $K$ whose support 
is $\ol C'$.
\par
Now $\CW = N_G(L)/L$ is isomorphic to the symmetric group $S_{n/d}$.
The irreducible character $E  = E_{\mu} \in S_{n/d}\wg$ corresponding to
$(C',\CE')$ under the generalized Springer correspondence 
is given by $\mu = (\la_1/d, \la_2/d, \dots)$.
\para{3.3.}
We fix an $F$-stable Borel subgroup $B$ of $G$ and an $F$-stable maximal
torus contained in $B$, where $B$ (resp. $T$) is the subgroup 
of $G$ consisting of upper triangular matrices  (resp. diagonal
matrices).  We fix $d$ as in 3.2, and put $t = n/d$.  
Let $P = LU_P$ be the parabolic subgroup of $G$ containing
$B$, where 
$L$ is the Levi subgroup of $P$ containing $T$  of type 
$A_{d-1}\times\cdots\times A_{d-1}$, ($t$-times).  Hence $P$, $L$ and
$U_P$ are all $F$-stable.  Let $(C, \CE)$ be the cuspidal pair in $L$ 
corresponding to $\r_0 \in A_L(v)\wg$ as in 3.2, and 
$(\ZC, \CL)$ the corresponding objects in $\Fl$.  The unipotent class
$C$ in $L$ can be identified with $C_1\times \cdots\times C_t$ in 
$SL_d\times\cdots\times SL_d$ with $C_i$ regular unipotent in
$SL_d$.  We choose $v = v_0 \in C^F$
so that $v_0$ is a product of split elements in $C^F_i$, and let 
$y_0 = v_0-1$ the corresponding element in $\ZC^F$ .  
Let $\wt A_L(v_0)$ be as in 1.3. 
Since $A_L(v_0)$ is abelian, $\r_0 \in A_L(v)\wg$ 
is linear.  We choose an extension $\wt\r_0$
so that $\wt\r_0(\s) = 1$.  This corresponds to an isomorphism  
$\vf_0: F^*\CE \isom \CE$ which induces the identity map on
the stalk $\CE_{v_0}$.
\par  
Let $\la = (\la_1, \dots, \la_r)$ be a
partition of $n$ such that all the $\la_i$ are divisible by $d$, and
$u = u_{\la}$  the split unipotent element in $G^F$.  As in the case 
of $(C,\CE)$, we choose an extension $\wt\r$ of $\r \in A_G(u)\wg$ 
corresponding to $\CE'$ by the condition that $\wt\r(\s) = 1$, and 
consider $\g = \g(v,\wt\r_0,u,\wt\r)$ as in 1.3.  Passing to 
the Lie algebra situation, we consider 
$y = y_{\la} \in  \Fg^F$ and $y_0 \in \ZC^F$.  Under this setting, 
we write $\g$ as $\g = \g(y_0, \wt\r_0, y, \wt\r)$.
We consider the subvariety $\CP_y$ of $G/P$ as given in (1.6.1).
As in 1.6, the map $\vf_0$ induces a linear isomorphism 
$\F$ on $H_c^{a_0+r}(\CP_y, \dot\CL)$.  We have
\begin{thm}  
Assume that $p$ is arbitrary, and let $G = SL_n$ with 
the standard Frobenius map $F$.
Then $\F$ acts on 
$H_c^{a_0+r}(\CP_y, \dot\CL) = H_c^{a_0+r}(\CP_y, \dot\CL)_{\r}$ as 
$q^{(a_0+r)/2}$ times identity.  In particular, we have 
$\g(y_0,\wt\r_0, y, \wt\r) = 1$.
\end{thm}
\para{3.5.}
The remainder of this section is devoted to the proof of Theorem 3.4.
Since the second statement easily follows from the first one, we
concentrate to the proof of the first statement. 
First we note that $\CP_y$ may be identified with the 
set $\CF_y$ of partial flags 
$$
D = (V_d \subset V_{2d} 
     \subset \cdots \subset V_{(t-1)d}),
$$
such that 
$D$ is $y$-stable and that $y$ induces a regular nilpotent transformation
on $V_{id}/V_{(i-1)d}$ for each $i \ge 1$.  (Here $V_j$ denotes a subspace
of $V$ with $\dim V_j = j$).   
\par
Let $\CG_y$ be the set of $d$-dimensional subspaces $V_d$ of $V$ such
that $V_d$ is $y$-stable and that $y$ acts as a regular nilpotent
transformation on $V_d$.  We have a natural surjective map 
$p : \CF_y \to \CG_y$ by $p(D) = V_d$.  Then $\CG_y$ is identified
with the variety $\BP(\Ker y^d) - \BP(\Ker y^{d-1})$; for each 
$v \in \Ker y^d -\Ker y^{d-1}$,
the space spanned by $v, yv, \dots y^{d-1}v$ gives an element in 
$\CG_y$.      
We have a filtration of $\CG_y$
$$
\CG_y = \CG_0 \supset \CG_1 \supset \cdots  , 
$$
where $\CG_i - \CG_{i+1} \simeq \BA^{s-i}$ with 
$\dim \CG_y = s = d(\dim\Ker y)-1$.
Here $\CG_i$ is defined by $\BP(U_i) - \BP(\Ker y^{d-1})$ for a
certain subspace $U_i$ of $\Ker y^d$ containing $\Ker y^{d-1}$ such
that $\Ker y^d = U_0 \supset U_1 \supset \cdots $.  
Let us   
choose a non-zero vector $w_i \in U_i - U_{i+1}$ for each $i$. We can
choose some $e_j$ as $w_i$.  As in the
case of $\CB_u$ for $GL_n$, one can define a map 
$f^{(i)}: \BA^{s-i} \to Z_{\wt G}(y), v \mapsto f^{(i)}_v$  
such that $f^{(i)}_v\cdot w_i = v$  
for $v \in U_i-U_{i-1}$, under the identification 
$\BP(U_i) - \BP(U_{i-1}) \simeq \BA^{s-i}$. (Here $\wt G$ denotes 
$GL_n$).  Let $V^{(i)}_d$ be the element in
$\CG_y$ corresponding to $w_i$.  Then $y$ induces a nilpotent
transformation $\ol y$ on $\ol V = V/V^{(i)}_d$, which corresponds to a
partition $\la'$ of $n-d$ obtained from $\la$  by replacing some 
$\la_j$ by $\la_j-d$.  Moreover,
$p\iv(V^{(i)}_d)$ is isomorphic to $\CF_{\ol y}$, the corresponding
variety for $SL(\ol V)$, under the correspondence 
$$
D = (V^{(i)}_d \subset  V_{2d} \subset \cdots \subset V_{(t-1)d}) \mapsto 
           \ol D = (\ol V_{2d} \subset \cdots \subset \ol V_{(t-1)d})  
$$
with $\ol V_{jd} = V_{jd}/V_d^{(i)}$.  As in the case of $GL_n$, by using 
the map $f^{(i)} : \AA^{s-i} \to Z_{\tilde G}(y)$, we have an isomorphism
\begin{equation*}
\tag{3.5.1}
p\iv(V_d^{(i)}) \times (\CG_i-\CG_{i+1}) \simeq p\iv(\CG_i - \CG_{i+1}),
\qquad (D, v) \mapsto f^{(i)}_v\cdot D.
\end{equation*}
Note that $\CF_y$ and $\CG_y$ have natural $\Fq$-structures inherited
from $G/P$.  Then $\CG_i$  are all $F$-stable, and the
isomorphism in (3.5.1) is $F$-equivariant.
\para{3.6.}
Let $Q$ be the maximal parabolic subgroup of $G$ containing $P$ of
type $A_{n-d-1}\times A_{d-1}$.  Let $\CG$ be the set of subspaces 
of dimension $d$ in $V$.
Then $\CG$ may be identified with $G/Q$ and $\CG_y$ is a locally
closed subvariety of $\CG$. The map $p: \CF_y \to \CG_y$ is obtained 
from the map $G/P \to G/Q$ by restricting it to $\CP_y$, which we also
denote by $p$. 
Now, $V_d^{(i)} \in \CG_y$ corresponds to   
$gQ \in G/Q$ for some $g = g_i \in G$ and 
$p\iv(V_d^{(i)})$ may be identified with $ p\iv(gQ)$, where
$$
p\iv(gQ) = \{ xP \in gQ/P \mid \Ad(x)\iv y \in \ZC + \Fn_P\}.
$$
We may choose $g$ so that $gP \in \CP_y$.
\par
We note that  $Q/P$ is isomorphic to $M/P_M$,
where $M$ is the subgroup of $G$ isomorphic to $SL_{n-d}$, 
and is isogeneous to a component of the Levi subgroup of $Q$ containing $T$.
Then $P_M = P \cap M$ is the parabolic
subgroup of $M$ of type $A_{d-1}\times \cdots \times A_{d-1}$, ($t-1$
factors), and $L_M = L \cap M$ is the Levi subgroup of $P_M$. 
The regular nilpotent orbit $\ZC$ in $\Fl$ can be written as 
$\ZC = \ZC_M \times \ZC_t$, where $\ZC_M$ is the regular nilpotent 
orbit in $\Lie L_M = \Fl_M$ 
and $\ZC_t$ is the regular nilpotent orbit in the $t$-th component of
$\Fl$.  Since $\Ad(g)\iv y \in \ZC + \Fn_P$, one can write 
$\Ad(g)\iv y = y' + z'$ with 
$y' \in \Fm$ and $z' \in \ZC_t + \Fn_Q$.  (Here $\Fm = \Lie M$ and 
$\Fn_{Q} = \Lie U_Q$). 
Set
\begin{align*}
\CP_{y'}^M &= \{ xP_M \in M/P_M \mid \Ad(x)\iv y' \in \ZC_M+ \Fn_{P_M}\}, \\  
\hat\CP_{y'}^M &= \{ x \in M \mid \Ad(x)\iv y'  \in \ZC_M + \Fn_{P_M}\}.
\end{align*}
We note that  
\par\medskip\noindent
(3.6.1) \ The map $xP_M \mapsto gxP$ gives an isomorphism 
$\CP_{y'}^M \simeq  p\iv(gQ)$.
\par\medskip 
In fact, since $M$ normalizes $C_tU_Q$, we have 
$\Ad(x)\iv z' \in \ZC_t + \Fn_Q$.
Then the condition $\Ad(gx)\iv y \in \ZC + \Fn_P$ is equivalent to 
the condition $\Ad(x)\iv y' \in \ZC_M + \Fn_{P_M}$.  (3.6.1) follows from this.
\par
By (3.6.1), one can define an injective map $\io: \CP_{y'}^M \to \CP_y$. 
Similarly,  
$\wh \CP_{y'}^M$ is isomorphic to
the set $\{ x' \in gM \mid \Ad(x')\iv y \in \ZC + \Fn_P\}$ which is a subset of 
$\wh \CP_y$.  Hence we have an injective map 
$\hat\io: \wh \CP_{y'}^M \to \wh \CP_y$.  Now it is easy to see that the 
following diagram commutes.
\begin{equation*}
\begin{CD}
\ZC @<\a<< \wh \CP_y @>\b>> \CP_y \\
@A\io' AA  @AA\hat\io A  @AA\io A \\
\ZC_M @<\a'<< \wh \CP_{y'}^M @>\b'>> \CP_{y'}^M.
\end{CD}
\tag{3.6.2}
\end{equation*}
Here the left vertical map is an injection 
$\io': \ZC_M \to \ZC, x \mapsto (x, y'')$, where $y''$ is the projection of 
$z' \in \ZC_t + \Fn_Q$ to $\ZC_t$, i.e., the projection of 
$\Ad(g)\iv y \in \ZC + \Fn_P$ 
on $\ZC_t$. The horizontal maps $\a', \b'$ are defined in a similar 
way as $\a$ and $\b$ by replacing $G$ by $M$. 
\par
Let $\CL$ and $\dot\CL$ be local systems on $\ZC$ and $\CP_y$,
respectively, as in 1.6.  We denote by $\CL_M$ and $\dot\CL_M$
similar objects for $\ZC_M$ and $\CP^M_{y'}$ as $\CL, \dot\CL$ for 
$\ZC$ and $\CP_y$.
Then $\CL_M$ coincides with $(\io')^*\CL$.  This implies, by (3.6.2),  
that 
\begin{equation*}
\tag{3.6.3}
\io^*\dot\CL = \dot\CL_M.
\end{equation*}
\par
Put
\begin{align*}
Y_y &= \{ xQ \in G/Q \mid \Ad(x)\iv y \in \Fm + \ZC_t + \Fn_Q \},  \\
\hat Y_y &= \{ x \in G \mid \Ad(x)\iv y \in \Fm + \ZC_t + \Fn_Q \}.
\end{align*}  
Then $Y_y$ is isomorphic to $\CG_y$.  
We consider the subset $Y_i$ of $Y_y$ corresponding to $\CG_i$.  Since 
$\CG_i - \CG_{i+1}$ coincides with the set 
$\{ f^{(i)}_v\cdot w_i \mid v \in \BA^{s-i} \}$,  
one can write as 
$Y_i-Y_{i+1} = \{ f^{(i)}_v g_iQ \mid v \in \BA^{s-i}\}$.
Then we have the following commutative diagram
\begin{equation*}
\begin{CD}
\ZC_t @<\wt\a << \wh Y_y @>\wt\b >> Y_y \\
@AAA    @AAA  @AAA         \\
\{ y''_i\}  @<<< f^{(i)}(\AA^{s-i})g_i @>>> Y_i-Y_{i+1}.
\end{CD}
\tag{3.6.4}
\end{equation*}
Here $y''_i = y''$ is as in (3.6.2), and
$\tilde\a(x)$ is the
$\ZC_t$-component of  $\Ad(x)\iv y \in \Fm + \ZC_t + \Fn_Q$, 
$\tilde\b(x) = xQ$. 
All the vertical maps are natural inclusions and the lower 
horizontal arrows are the restrictions of upper ones.  Note that 
the right lower horizontal map is an isomorphism since 
$Y_i - Y_{i+1} \simeq \BA^{s-i}$.
\par
Let $\CL_t$ be the cuspidal local system on $\ZC_t$.  Then we have a
local system $\dot\CL_t$ on $Y_y$ by the condition that 
$\wt\a^*\CL_t = \wt\b^*\dot\CL_t$.  Since $\CL_t$ is a local
system of rank 1, it follows from (3.6.4) that
\par\medskip\noindent
(3.6.5) \  The restriction of $\dot\CL_t$ on $Y_i-Y_{i+1}$ is the constant
sheaf $\Ql$.
\par\medskip 
We now consider the commutative diagram
\begin{equation*}
\begin{CD}
\ZC @<\a << \wh \CP_y @>\b >> \CP_y \\
@AAA   @AAA   @AAA   \\
\ZC_M\times \{y_i''\}  @<\a''<< \b\iv(p\iv(Y_i-Y_{i+1})) 
        @>\b''>> p\iv(Y_i-Y_{i+1}).
\end{CD}
\tag{3.6.6}
\end{equation*}     
Here all the vertical maps are natural inclusions, and the horizontal
maps $\a''$ and $\b''$ are the restrictions of $\a$ and $\b$.
By (3.5.1), we have 
\begin{align*}
\tag{3.6.7}
p\iv(Y_i-Y_{i+1}) &\simeq \CP_{y'}^M \times (Y_i-Y_{i+1}), \\ 
       \b\iv(p\iv(Y_i-Y_{i+1})) 
         &\simeq \wh \CP_{y'}^M \times f^{(i)}(\AA^{s-i})g_i, 
\end{align*}
and under the above isomorphisms, the maps $\a'', \b''$ are given as 
$$
\a''(x, f^{(i)}_vg_i) = (\a'(x), y_i''), 
              \qquad \b''(x, f^{(i)}_vg_i)) = (\b'(x), v) 
$$
for $x \in \CP^M_{y'}, v \in Y_i-Y_{i+1} \simeq \BA^{s-i}$.
\par
Now the restriction of $\CL$ to $C_M\times \{y''_i\}$ is 
a local system $\CL_M\boxtimes \Ql$.  Hence by
making use of (3.6.6) and (3.6.7), we have
\par\medskip\noindent
(3.6.8) \ Under the isomorphism 
$p\iv(Y_i-Y_{i+1}) \simeq \CP_{y'}^M \times (Y_i-Y_{i+1})$,  
the restriction of 
$\dot\CL$ on $p\iv(Y_i-Y_{i+1})$ coincides with 
$\dot\CL_M \boxtimes \Ql$.
\par\medskip
It follows from (3.6.8) that we have an isomorphism
\begin{equation*}
\tag{3.6.9}
H_c^k(p\iv(Y_i-Y_{i+1}), \dot\CL) \simeq 
                      H_c^{k'}(\CP_{y'}^M, \dot\CL_M),
\end{equation*}
where $k \equiv k'\pmod 2$.  Then using the locally trivial filtration 
of $\CP_y = p\iv(Y_0) \supset p\iv(Y_1) \supset \cdots$, and by induction 
on the rank of $G$, we see that
\begin{equation*}
\tag{3.6.10}
H_c^{\odd}(p\iv(Y_i),\dot\CL) = 0
\end{equation*}
for any $i \ge 0$.
\para{3.7.}
We are now ready to prove Theorem 3.4.
Put $m = a_0 + r$.  First we note the following.
\par\medskip\noindent
(3.7.1) \  $H_c^m(\CP_y, \dot\CL) = H_c^m(\CP_y, \dot\CL)_{\r}$, 
and the map $\F$ acts on $H_c^m(\CP_y, \dot\CL)$ as a
scalar multiplication.  
\par\medskip
In fact, it follows from section 2 that  $H_c^m(\CP_y,\dot\CL)$ has 
a natural structure of $\CW \times A_G(y)$-module, which is compatible
with the isomorphisms (1.4.3) and (1.6.2).  Hence by the generalized 
Springer correspondence, it is decomposed as 
\begin{equation*}
H_c^m(\CP_y,\dot\CL) \simeq \bigoplus_{\r' \in A_G(y)\wg}
                        V_{y,\r}\otimes\r',
\end{equation*}
where $V_{y,\r'}$ is an irreducible $\CW$-module whenever 
it is non-zero.
Now the explicit description of the generalized Springer correspondence
in the case of $SL_n$ (see. 3.2) shows that
$\r$ is the unique character such that $V_{y,\r} \ne 0$.  Hence 
$H_c^m(\CP_y, \dot\CL) = H_c^m(\CP_y, \dot\CL)_{\r}$.  Since $A_G(y)$
is abelian, $H_c^m(\CP_y, \dot\CL)$ is an irreducible $\CW$-module. 
It is easy to see that the map $\F$ on $H^m_c(\CP_y,\dot\CL)$ 
commutes with the action of $\CW$.  Hence $\F$ is a scalar
multiplication,  and so (3.7.1) holds.
\par\medskip
Note that in the discussion of 3.5 and 3.6, $\CG_y, Y_y,$ etc.  have 
natural $\Fq$-structures.  We may choose the filtration of $\CG_y$ 
and $Y_y$ compatible with the $\Fq$-structure, i.e., all the $Y_i$
and $\wh Y_y$ are $F$-stable.  Then all the diagrams and formulas
there hold with $\Fq$-structure.     
We consider the top piece $Y_0-Y_1$ of the filtration of 
$Y_y$.  In this case, we may choose $g = g_0 = 1$ in the discussion in 3.6, 
and so $y$ is decomposed as $y = y' + z'$ with $y' \in \Fm$ and 
$z' \in \ZC_t + \Fn_Q$.  Hence $y'$ (resp. $y''$ ) is the projection of 
$y$ on $\Fm$ (resp. on $\ZC_t$).  Since $y$ is a split element, $y', y''$
are also split.  Let 
\begin{equation*}
m' = (\dim M - \dim \ZC_{y'}) - (\dim L_M - \dim \ZC_M),
\end{equation*}
where $\ZC_{y'}$ is
the nilpotent orbit in $\Fm$ containing $y'$.  Since $\ZC$ is the
regular nilpotent orbit in $\Fl$, we see easily that $m = 2 \dim \CB_y$, where 
$\CB_y$ is the variety of Borel subgroups whose Lie algebra contains
$y$.  Similarly, 
we have $m' = 2 \dim\CB_{y'}^M$.  Then by using the locally trivial 
filtration of $\CB_y$, we see  
that 
\begin{equation*}
\tag{3.7.2}
m - m' = 2d\dim\Ker y = 2s.
\end{equation*}
In fact, assume that $y = y_{\la}$ with 
$\la = (\la_1, \dots, \la_k)$.  By using the locally trivial
filtration arising from the maximal parabolic subgroup $P_1$ of $G$ with
Levi subgroup $L_1$ of type $A_{n-2}$, one
obtains that $\dim \CB_y - \dim \CB^{L_1}_{y_1} = \dim\Ker y$, where
$y_1$ is a nilpotent element in $\Lie L_1$ of type 
$\la' = (\la_1, \dots, \la_k-1)$.
In the same way, one can find a 
nilpotent element $y_2 \in \Lie L_2$ with type
$(\la_1, \dots, \la_k-2)$ such that 
$\dim \CB^{L_1}_{y_1} - \dim \CB^{L_2}_{y_2} = \dim\Ker y$, where
$L_2$ is a Levi subgroup of the maximal parabolic subgroup $P_2$ of $L_1$.
Repeating this procedure, one can find similar formulas for 
$L_1 \supset L_2 \supset \cdots \supset L_d$ with  
$\CB^{L_d}_{y_d} = \CB^M_{y'}$.  (3.7.2) follows from this.
\par
 Since $Y_0-Y_1 \simeq \BA^s$, we
have an isomorphism with $\Fq$-structures
\begin{equation*}
\tag{3.7.3}
H_c^m(p\iv(Y_0-Y_1), \dot\CL) \simeq 
             H_c^{m'}(\CP_{y'}^M, \dot\CL_M)[s]
\end{equation*}
as a special case of (3.6.9), 
where $[s]$ is the Tate twist.  (The compatibility of the
Frobenius actions comes from the discussion in 3.6 by noticing that
$y''$ is a split element in $\ZC_t$.)  Let $\F_M$ be the map on 
$H_c^{m'}(\CP_{y'}^M, \dot\CL_M)$ defined in a similar way as $\F$.
By induction on the rank of
$G$, we may assume that $\F_M$ acts on $H_c^{m'}(\CP_{y'}^M, \dot\CL_M)$ 
as a scalar multiplication by $q^{m'/2}$.  Then by (3.7.3), $\F$ acts
on $H_c^m(p\iv(Y_0-Y_1),\dot\CL)$ as a scalar multiplication by 
$q^{m/2}$.  Now by using the cohomology long exact sequence with respect
to the closed immersion $p\iv(Y_1) \subset p\iv(Y_0) = \CP_u$,
together with (3.6.10), we see that the natural map
$$
H_c^m(p\iv(Y_0-Y_1), \dot\CL) \longrightarrow H_c^m(\CP_y,\dot\CL)
$$
is injective.  This proves the theorem since $\F$ acts on 
$H_c^m(\CP_y,\dot\CL)$ by a scalar multiplication by (3.7.1).

\par\medskip
\section{$G = SL_n$ with $F$ of non-split type}  
\para{4.1.}
In this section, we assume that $G = SL_n$ is as in section 3,
and that $p$ is large enough so that the argument in section 2
can be applied (e.g., $p > 3(n-1)$).
Let $F = \s F_0$ be the twisted Frobenius map on $G$, where
$F_0$ is the standard Frobenius map over $\Fq$ as in 3.1, and $\s$ is the 
graph automorphism on $G$ of order 2.  Here we take $\s: G \to G$
defined by $\s(g) = w_0{}^tg\iv w_0\iv$ for $g \in G$
($w_0$ is the permutation matrix in $GL_n$ corresponding to 
the longest element in $S_n$, and ${}^tg$ means the
transpose of the matrix $g = (g_{ij})$).  Then $B$ and $T$ in 3.3
are $F$ and $F_0$-stable.
\par
Unipotent classes in $G$ are all $F$-stable.  In order to describe
elements in $C^F$ for each unipotent class $C$, we introduce
a sesqui-linear form as follows.  Let $V \simeq k^n$ be as in 3.1, and
$V_0$ the $\BF_{q^2}$-subspace of $V$ generated by $\{ e_i\}$.  
We define a sesqui-linear form $\lp \ ,\ \rp$ on $V_0$ by
$\lp \sum_i a_ie_i, \sum_j b_je_j\rp = \sum_i a_ib^q_{n-i}$.  Then it
is easy to see that for $g \in G^{F_0^2}$, $g \in G^F$ if and only if
$\lp gv, gw\rp = \lp v, w\rp$ for any $v,w \in V_0$.  Let $\Fg$ be  
the Lie algebra of $G$, on which $F$ acts naturally.  Then for 
$x \in \Fg^{F^2_0}$, $x \in \Fg^F$ if and only if 
$\lp xv, w\rp + \lp v, xw\rp = 0$ for any $v,w \in V_0$.
\para{4.2.}
For a partition $\la = (\la_1, \dots, \la_r)$ of $n$, we shall 
construct a nilpotent element $y_{\la} \in \Fg^F$.  First we note that
there exist basis vectors 
\begin{equation*}
f^{(i)}_j \quad (1 \le i \le r, 1 \le j \le \la_i)
\end{equation*}
of $V_0$ satisfying the property that 
\begin{equation*}
\lp f^{(i)}_j, f^{(i')}_{k}\rp = \begin{cases}
        1     &\quad\text{ if $i = i',j+k = \la_i+1$ and $j \ne k$}, \\
        \pm 1 &\quad\text{ if $i = i', j+k = \la_i+1$ and $j = k$}, \\      
        0     &\quad\text{ otherwise}.
                                 \end{cases}
\end{equation*} 
In fact, we can choose $f^{(i)}_j = e_k$ for some
$k$ if $\la_i$ is even.  If $\la_i$ is odd, put $\la_i = 2t_i+1$.  
Then $f^{(i)}_j$ is of the form $e_k$ if $j \ne t_i+1$, and 
we can choose $f^{(i)}_{t_i+1}$ from one of the vectors 
$e_k \pm \frac 1 2 e_{n-k+1}$ with $2k \ne n+1$ 
(note that $p > 2$), or $e_l$ with $n = 2l + 1$. 
\par
Put $t_i = [\la_i/2]$ ($[\ \ ]$ is the Gauss symbol) for each $\la_i$.
We now define a nilpotent transformation $y_{\la} \in \Fg^{F^2}$ 
on $V_0$ by
\begin{equation*}
y_{\la}f^{(i)}_j = \begin{cases}
        f^{(i)}_{j+1} &\quad\text{ if } 1 \le j \le t_i-1, \\
        \ve_if^{(i)}_{j+1} &\quad\text{ if } j = t_i, \\
        -f^{(i)}_{j+1} &\quad\text{ if } t_i+1 \le j \le \la_i-1, \\
         0                  &\quad\text{ if } j = \la_i,
                \end{cases}
\end{equation*}
where $\ve_i = 1$ if $\la_i$ is even, and 
$\ve_i = \lp f^{(i)}_{j+1}, f^{(i)}_{j+1}\rp$ if $\la_i$ is odd. 
Then 
\begin{equation*}
\tag{4.2.1}
\{ y_{\la}^jf^{(i)}_1 \mid 1 \le i \le r, 0 \le j \le \la_i-1\} 
\end{equation*}
gives a basis of $V_0$ satisfying the relation
\begin{align*}
\tag{4.2.2}
\lp y_{\la}^jf^{(i)}_1, y_{\la}^{\la_i-j+1}f^{(i)}_1\rp = (-1)^ja_i, \quad
(a_i = \pm 1)
\end{align*}
and $\lp y_{\la}^jf^{(i)}_1, y_{\la}^kf^{(i')}_1\rp = 0$ for all other
pairs.  It follows from this that $y_{\la} \in \Fg^F$.
\par
Let $d$ be as in 3.2, and assume that $d \ge 2$.
We assume that the partition $\la$ satisfies the condition that 
all the parts $\la_i$ are divisible by $d$.  
We shall construct a nilpotent element $y_1 \in \Fg^F$ of type
$\n = (d, \dots, d)$ associated to $y_{\la}$.  We define a 
map $y_1$ on $V_0$ by 
\begin{equation*}
\tag{4.2.3}
y_1f^{(i)}_j = \begin{cases}
      0    &\quad\text{ if } j \equiv 0 \pmod d, \\
      y_{\la}f^{(i)}_j &\quad\text{ otherwise}.
                \end{cases}   
\end{equation*}
Then in view of (4.2.2), it is easy to check that $y_1$ leaves
the form $\lp\ ,\ \rp$ invariant, and we have $y_1 \in \Fg^F$.  
\para{4.3.} 
Let $L$ be a Levi subgroup of the standard parabolic subgroup
$P$ of $G$ of type $A_{d-1}\times \cdots\times A_{d-1}$ 
($t = n/d$-factors).  
(Here $P$ and $L$ are as in 3.1 with respect 
to $F_0$, $P$ is $\s$-stable, $\s$ permutes the $i$-th factor and
$(t-i+1)$-th factor, etc. )  Thus $P$ and $L$ are $F$-stable.  
Let $\ZC$ be the regular 
nilpotent orbit in $\Fl$.  We choose a representative $y_0 \in \ZC^F$
in the following way;  we define
a basis $\{ e^{(i)}_j \mid 1 \le i \le t, 1 \le j \le d\}$ of 
$V_0$ by
$e^{(i)}_j = e_{(i-1)d + j}$.
Then in the case where $t$ is even, or $t$ is odd and  
$i \ne (t+1)/2$, we define
\begin{equation*}
y_0e^{(i)}_j = \begin{cases}
      e^{(i)}_{j+1}  &\quad\text{ if } 1 \le i \le [t/2], j \ne d,\\
      -e^{(i)}_{j+1} &\quad\text{ if } 
             t - [t/2]+1 \le i \le t, j \ne d, \\
         0           &\quad\text{ if }  j = d  
               \end{cases}
\end{equation*}
If $t$ is odd and $i = (t+1)/2$, let $V_1$ be the subspace 
of $V$ spanned by $e^{(i)}_j$ with $1 \le j \le d$.  
We define $y_0|_{V_1}$ as a regular nilpotent element 
$y_{\la} \in \Fs\Fl_d^F$ as in 4.2.  
\par
Let $(\ZC, \CL)$ be the cuspidal pair in  $\Fl$ corresponding to
an $F$-stable character $\r_0$ of $A_L(y_0)$.  
We have a natural homomorphism $A_L(y_0) \to A_G(y_0)$.  Since
$A_G(y_0)$ is a cyclic group of order $d$, this gives 
an isomorphism compatible with $F$-action.
Thus $\r_0$ is regarded as an $F$-stable character of $A_G(y_0)$.
Since $y_0$ and $y_1$ are conjugate under $G$, there exists 
$c_1 \in A_G(y_0)$ (up to $F$-conjugacy) such that $y_1$ is obtained 
from $y_0$ by twisting by $c_1$. 
Since $\r_0$ is $F$-stable, the value $\r_0(c_1)$ is well-defined.  This 
value is determined by $y_1$, hence by $y_{\la}$, which we denote
by $\e_{\la}$.  
Let $\g(y_0, \wt\r_0, y_{\la}, \wt\r)$ be the scalar defined
by choosing the extensions $\wt\r_0, \wt\r$ in a similar way as 
the case of split $F$ (cf. 3.3). 
Put $m = a_0 + r$ as before.   
We have the following theorem.
\begin{thm} 
Assume that $p$ is large enough so that Dynkin-Kostant theory can be
applied.  
Let $w_0$ be the longest element in $\CW$.  Then $\F w_0$ acts on
$H^{m}_c(\CP_{y_{\la}}, \dot\CL) = 
           H^{m}_c(\CP_{y_{\la}}, \dot\CL)_{\r}$ 
as a scalar multiplication
by $\e_{\la}(-q)^{m/2}$.  In particular, 
$\g(y_0, \wt\r_0, y_{\la}, \wt\r) = \e_{\la}(-1)^{m/2}$.
\end{thm} 
\para{4.5.}
The remainder of this section is devoted to the proof of the theorem.
If we notice that the preferred extension $\wt V_E$ of $V_E$ is 
given by defining the action of $\s \in \wt\CW$ 
by the action of $w_0 \in \CW$, 
the second statement follows easily from the first one.  So 
we concentrate the proof of the first statement.
For $y_1$ of type $(d, \dots, d)$, we construct a $\Fs\Fl_2$-triple
$\{ y_1, y_1^-, h_1\}$ as follows.  On each Jordan block, $y_1$ can be
expressed as a matrix of degree $d$ with respect to the basis 
in (4.2.1) as 
\begin{equation*}
Y = \begin{pmatrix}
       0 &   &  &   \\
       1 &  \ddots  &  &   \\
         &  \ddots   & 0 &  \\
         &        & 1 & 0   \\
          \end{pmatrix}.
\end{equation*} 
We define matrices $Y^-, H$ of degree $d$ by
\begin{align*}
Y^- &= \begin{pmatrix}
         0 & 1\cdot(d-1) &         &         & \\
          &  0      & 2(d-2) &         &  \\
          &         &   0     &  \ddots &   \\ 
          &         &         &  \ddots & (d-1)\cdot 1  \\
          &         &         &         &  0
      \end{pmatrix},   \\
\end{align*}
\begin{align*}
H &= \begin{pmatrix}
          1-d &     &        &    \\
              & 3-d &        &    \\
              &     & \ddots &  \\
              &     &        & d-1
     \end{pmatrix}.
\end{align*}
Then $[H, Y] = 2Y, [H, Y^-] = -2Y, [Y, Y^-] = H$.  
Thus by combining these matrices for all the Jordan blocks,
one obtains $y_1^-, h_1 \in \Fg$ satisfying the property that 
$[h_1,y_1] = 2y_1, [h_1, y_1^-] = -2y_1^-, [y_1, y_1^-] = h_1$
as asserted.  It follows from the construction, we see easily that
$y_1^-, h_1 \in \Fg^F$. 
\par 
We define a transversal slice $\vS$ with respect to the orbit through  
$y_1$ in $\Fg$ by $\vS = y_1 + Z_{\Fg}(y_1^-)$.  Hence $\vS$ is
$F$-stable.  We have the following lemma.
\begin{lem}
Let $y_{\la}$ be as in 4.2.  Then we have $y_{\la} \in \vS^F$.
\end{lem}
\begin{proof}
We write $y_{\la}$ as $y_{\la} = y_1 + y$.  It is enough to show
that  $y \in Z_{\Fg}(y_1^-)$.
Now $y$ is a nilpotent transformation on $V_0$ determined by the
condition that 
\begin{equation*}
y: y_{\la}^jf^{(i)}_1 \mapsto y_{\la}^{j+1}f^{(i)}_1
\end{equation*} 
for $j \equiv 0 \pmod d$, and it maps all other $y_{\la}^jf^{(i)}_1$ to 0.
Since $y_1^-$ maps $y_{\la}^jf^{(i)}_1$ to $y_{\la}^{j-1}f^{(i)}_1$
up to scalar if $j\not\equiv 1 \pmod d$, and to 0 if $j \equiv 1\pmod d$.
 we see easily check 
that $y_1^-\circ y = y\circ y_1^- = 0$ on $V_0$.  
Hence $y \in Z_{\Fg}(y_1^-)$. 
\end{proof} 
\para{4.7.}
By using the $\Fs\Fl_2$-triple $\{ y_1, y_1^-, h_1\}$, one can define 
a Lie algebra homomorphism $\f_1: \Fs\Fl_2 \to \Fg$ as in 2.7.
The construction of $\Fs\Fl_2$-triple given in 4.5 works well for $y_{\la}$
in general, and one gets the $\Fs\Fl_2$-triple containing $y_{\la}$. 
We denote by $\f_{\la}$ the homomorphism $\Fs\Fl_2 \to \Fg$ obtained 
from it.  Thus $Z_G(\f), M_G(\f)$, are defined as in 2.7 for 
$\f = \f_1, \f_{\la}$.  
\par
Let $\pi: \dot\Fg \to \Fg$ be as in 2.8.  Then $\CP_y \subset \dot\Fg$, 
and the local system $\dot\CL$ on $\CP_y$ can be extended to a local system on
$\dot\Fg$ (cf. 2.8), which we denote also by $\dot\CL$ as
in 2.8.   Put $K_1 = \pi_!\dot\CL$.  $K_1$ is essentially the same as 
$K = \pi_!\dot\CL^*$ in 2.8, and so $K_1[\d]$ is a perverse sheaf on 
$\Fg$ with 
a canonical $\CW$-action, where $\d$ is as in 2.8.   
By making use of the transversal slice $\vS$, we show the following
proposition.
\begin{prop}  
There exist natural maps of $\CW$-modules, which make the
following diagram commutative.
\begin{figure}[h]
\setlength{\unitlength}{1mm}
\begin{picture}(80, 40)(5,-5)
\put(20, 20){\makebox(20,20){$\HH^i(\Fg, K_1)$}}
\put(44, 30){\vector (1,0){10}}
\put(44, 27){\makebox(10,10) {\footnotesize{$\pi_1$}}}
\put(60, 20){\makebox(20,20){$ H_c^i(\CP_{y_1}, \dot\CL)$}}
\put(30, 25){\vector (0,-1){10}}
\put(20, -1){\makebox(20,20){$H_c^i(\CP_{y_{\la}}, \dot\CL)$}}
\put(58, 25){\vector (-3,-2) {15}}
\put(22, 15){\makebox(10,10){\footnotesize{$\pi_{\la}$}}}
\put(52, 15){\makebox(10,10){\footnotesize{$\xi_{\la}$}}}
\end{picture}
\end{figure}
\par
\vspace{-1cm}
\noindent
Moreover, the map $\xi_{\la}$ is equivariant with respect to the
actions of $\F$ on both cohomologies.
\end{prop}
\newpage
\begin{proof}
By the inclusion $\{ y_{\la}\} \subset \vS \subset \Fg$, we have
the canonical maps
\par\medskip
\begin{figure}[h]
\setlength{\unitlength}{1mm}
\begin{picture}(80, 40)(5,-5)
\put(20, 20){\makebox(20,20){$\HH^i(\Fg, K_1)$}}
\put(44, 30){\vector (1,0){10}}
\put(60, 20){\makebox(20,20){$ \HH^i(\vS, K_1)$}}
\put(30, 25){\vector (0,-1){10}}
\put(20, -1){\makebox(20,20){$\CH^i_{y_{\la}}(K_1)$}}
\put(58, 25){\vector (-3,-2) {15}}
\end{picture}
\end{figure}
\vspace{-32mm}
\noindent
(4.8.1)
\par
\vspace{20mm}
Since $K_1$ is a $\CW$-complex with respect to the trivial action of 
$\CW$ on $\Fg$, the above maps are $\CW$-equivariant.
Since $K_1= \pi_!\dot\CL$, we have 
\begin{equation*}
\CH^i_{y_{\la}}(K_1) \simeq H^i_c(\CP_{y_{\la}}, \dot\CL)
\end{equation*}
by the proper base change theorem.
On the other hand, $K_1[\d]$ is a perverse sheaf on $\Fg$. 
Since the morphism $G \times \vS \to \Fg$ is 
smooth with all fibres of pure dimension equal to $\dim Z_G(y_1)$,
by a similar argument as in [L6, 3.2], $K_1[\dim \vS]|_{\vS}$ is
a perverse sheaf on $\vS$.  
$\vS$ is stable under the $\BG_m$-action ($t: x \mapsto t^{i-2}x$ for
each $x \in \Fg_i$ with respect to the grading $\Fg = \bigoplus \Fg_i$
associated to $\f_1:\Fs\Fl_2 \to \Fg$), and contracts to $y_1 \in
\vS$.  Since $K_1$ is $\BG_m$-equivariant, the canonical map
$\HH^i(\vS, K_1) \to \CH^i_{y_1}(K_1)$ gives rise to an isomorphism
\begin{equation*}
\HH^i(\vS, K_1) \simeq \CH^i_{y_1}(K_1) \simeq H^i_c(\CP_{y_1}, \dot\CL).
\end{equation*} 
The proposition follows from this.
\end{proof}
\par
For the special case where $i = 0$, we have the following
more precise result.
\begin{lem} 
The maps $\pi_1, \pi_{\la}$ in Proposition 4.8 give 
isomorphisms
\begin{equation*}
\HH^0(\Fg, K_1) \simeq H^0_c(\CP_{y}, \dot\CL) \simeq \vG(\ZC,\CL) 
\end{equation*}
for $y = y_1, y_{\la}$.  
In particular, 
$\xi_{\la}^0: H_c^0(\CP_{y_1}, \dot\CL) \to H_c^0(\CP_{y_{\la}},\dot\CL)$
is an isomorphism.
\end{lem}
\begin{proof}
We consider the following commutative diagram
\begin{equation*}
\begin{CD}
\ol\ZC @<\ol\a << \wh{\dot\Fg'} @>\ol\b>> \dot\Fg' @>\ol{\pi} >> \Fg \\
@A \bar j AA                @A \wh j AA  @AA j A     @AA\id A \\
\ZC @<\a <<    \wh{\dot \Fg} @>\b >>  \dot\Fg  @> \pi >> \Fg
\end{CD}
\tag{4.9.1}
\end{equation*}
where the lower horizontal maps are as in 2.8, and
\begin{align*}
\dot\Fg' &= \{ (x, gP) \in \Fg\times G/P \mid 
          \Ad(g\iv)x \in \ol\ZC + \Fz + \Fn_P\},   \\
\wh{\dot\Fg'} &= \{ (x,g) \in \Fg\times G \mid  
               \Ad(g\iv)x \in \ol\ZC + \Fz + \Fn_P\} 
\end{align*}
and $\ol\a, \ol\b$ and $\ol{\pi}$ are maps defined in a similar
way as $\a,\b$ and $\pi$.  
$\wh j, j$ are open immersions, and $\ol{\pi}$ is proper.
Now the local system $\dot\CL$ on $\dot\Fg$ is determined by the
condition that $\a^*\CL = \b^*\dot\CL$. 
Since the square in the left hand side in (4.9.1) is cartesian,
$\wh j_!(\a^*\CL) \simeq \ol\a^*(\bar j_!\CL)$.
The middle square is also cartesian, and we have 
\begin{equation*}
\ol\b^*(j_!\dot\CL) \simeq \wh j_!(\b^*\dot\CL) \simeq 
     \wh j_!(\a^*\CL) \simeq \ol\a^*(j_!\CL).
\end{equation*}
By the definition of the direct image with compact support, 
we have $\pi_!\dot\CL = \ol{\pi}_*(j_!\dot\CL)$.
Then
\begin{align*}
\HH^0(\Fg, K_1) &= \HH^0(\Fg, \pi_!\dot\CL) \\
              &= \HH^0(\Fg, \ol{\pi}_*(j_!\dot\CL)) \\
              &\simeq H^0(\dot\Fg', j_!\dot\CL). \\
\end{align*} 
Similarly, by using the open immersion 
\begin{equation*}
j: \CP_y \hra \ol\CP_y 
    = \{ gP \in G/P \mid \Ad(g\iv)y \in \ol\ZC + \Fn_P \}, 
\end{equation*}
we see that 
$H^0_c(\CP_y, \dot\CL) \simeq H^0(\ol\CP_y, j_!\dot\CL)$
for $y = y_1, y_{\la}$.  It follows that the maps $\pi_1, \pi_{\la}$ in 
Proposition 4.8 for $i = 0$ are nothing but the restriction 
$\vG(\dot\Fg', j_!\dot\CL) \to \vG(\ol\CP_y, j_!\dot\CL)$
of the global section of the sheaf $j_!\dot\CL$ on $\dot\Fg'$ 
for $y = y_1, y_{\la}$.
But since $\ol\b^*(j_!\dot\CL) \simeq \ol\a^*(\bar j_!\CL)$, we have
\begin{equation*}
\vG(\dot\Fg', j_!\dot\CL) \simeq \vG(\ol\ZC, \bar j_!\CL) 
    \simeq \vG(\ol{\CP}_y, j_!\dot\CL). 
\end{equation*}
Finally, we note that $\bar j_!\CL \simeq \bar j_*\CL$ since
$\CL$ is the cuspidal local system and so is clean ([L7, 2.2]).  
Hence 
\begin{equation*}
\vG(\ol\ZC, \ol j_!\CL) \simeq \vG(\ol\ZC, \ol j_*\CL) 
   \simeq \vG(\ZC, \CL)
\end{equation*} 
as asserted.
\end{proof}
\para{4.10.}
Let $\f_0: \Fs\Fl_2 \to \Fl \subset \Fg$ be the Lie algebra 
homomorphism such that $\f_0\binom{01}{00} = y_0$ constructed as in 
4.3.  Thus $\f_0$ is $F$-equivariant with respect to the twisted 
Frobenius action on $\Fs\Fl_2$.
Put $G_0 = Z_G^0(\f_0)$ and $T_0 = Z_L^0(\f_0)$.  Then $G_0$ and
$T_0$ are $F$-stable. 
It is checked that $G_0$ is isomorphic to $SL_t$, 
and $F$ acts as a twisted Frobenius endomorphism on $SL_t$. 
By (2.9.1) we have $T_0 \simeq Z_L^0$, and under the identification 
$G_0 \simeq SL_t$,
$T_0$ coincides with a  maximally split maximal
torus of $SL_t$, and $\CW = N_G(Z_L^0)/L$ is naturally isomorphic 
to the Weyl group of $G_0$ with respect
to $T_0$.  
\par
$F$ acts naturally on $\CW \simeq S_t$, as a conjugation 
by $w_0 \in \CW$, where $w_0$ is the longest element in $\CW$.  
Thus $Fw_0$ acts trivially on $\CW$.   
By 2.17, $F$ acts naturally on 
$H^*_{T_0} = \bigoplus_iH^{2i}_{T_0} \simeq S(\Fh^*)$, 
where $\Fh^* = \Ql\otimes_{\BZ}X(T_0)$.
$\CW$ also acts on $H^*_{T_0}$, which coincides with the  
action of $\CW$ on $S(\Fh^*)$ induced from the action of $\CW$ on 
$X(T_0)$ (cf. (2.3.2)).
We have the following lemma.
\begin{lem} 
$Fw_0$ acts on $H^{2i}_{T_0}$ as a scalar multiplication by $(-q)^i$.
\end{lem}
\begin{proof}
$Fw_0$ commutes with the graded algebra structure of $H^*_{T_0}$.
Since $H^*_{T_0}$ is generated by $H^2_{T_0}$, it is enough to show
that $Fw_0$ acts on $H^2_{T_0}$ as a scalar multiplication by $-q$. 
We show this by modifying the arguments used in the proof of
Lemma 2.4.  Let $\vG$ be as in the proof of Lemma 2.4 (with respect to
$G_0$).  We consider the locally trivial fibration   
$f: T_0\bs\vG \to G_0\bs\vG$.  We may assume that $\vG$ is defined 
over $\Fq$, and $f$ is $F$-equivariant.   We consider the
spectral sequence 
\begin{equation*}
\tag{4.11.1}
H^p(G_0\bs\vG, R^qf_*\Ql) \Rightarrow H^{p+q}(T_0\bs\vG), 
\end{equation*}
which have natural actions of $\CW$ (cf. 2.4) and of $F$.
Let $\th$ be the reflection representation of $\CW$.  Then (4.11.1)
implies a spectral sequence 
\begin{equation*}
H^p(G_0\bs\vG, R^qf_*\Ql)_{\th} \Rightarrow H^{p+q}(T_0\bs\vG)_{\th},
\end{equation*}
where $X_{\th}$ denotes the $\th$-isotypic subspace for a $\CW$-module $X$.
As in 2.4, we have 
\begin{equation*}
H^p(G_0\bs\vG, R^qf_*\Ql) \simeq 
        H^p(G_0\bs\vG)\otimes H^q(T_0\bs G_0),
\end{equation*}
and so
\begin{equation*}
H^p(G_0\bs\vG, R^qf_*\Ql)_{\th} \simeq 
        H^p(G_0\bs\vG)\otimes H^q(T_0\bs G_0)_{\th}
\end{equation*}
since $\CW$ acts trivially on $H^p(G_0\bs\vG)$.
Now it is known that 
$\bigoplus_iH^{2i}(T_0\bs G_0)$ is a graded regular representation of 
$\CW$, and that 
\begin{equation*}
H^q(T_0\bs G_0)_{\th} = \begin{cases}
                      H^q(T_0\bs G_0)  &\quad\text{ if } q = 2, \\
                      0                &\quad\text{ if } q < 2.
                  \end{cases} 
\end{equation*}
Since $H^*(T_0\bs \vG) = H^*_{T_0} \simeq S(\Fh^*)$, we have 
$H^2(T_0\bs \vG)_{\th} = H^2(T_0\bs \vG)$.  
Moreover, $H^0(G_0\backslash\vG) = H^0_G = \Ql$ be Lemma 2.4.
It follows that
\begin{equation*}
H^2(T_0\bs \vG) \simeq H^2(T_0\bs G_0).
\end{equation*} 
This isomorphism is compatible with the actions of $F$ and $\CW$.
It is well-known that $Fw_0$ acts as a scalar multiplication $-q$
on $H^2(T_0\bs G_0) = H^2(B_0\bs G_0)$, where $B_0$ is the $F$-stable
Borel subgroup of $G_0$ containing $T_0$. 
Hence $Fw_0$ acts similarly on $H^2_{T_0} = H^2(T_0\bs\vG)$.
This proves the lemma.
\end{proof}
\para{4.12.}
We consider the equivariant homology 
$H^{M^0(y_{\la})}_*(\CP_{y_{\la}}, \dot\CL^*)$, where
$M^0(y_{\la}) = M_G^0(y_{\la})$.  By results in 
Section 2, the graded Hecke algebra $\BH = \BS\otimes \Ql[\CW]$ 
acts on $H^{M^0(y_{\la})}_*(\CP_{y_{\la}},\dot\CL^*)$, where
$\BS = S(\Fh^*)\otimes\Ql[\Br]$ as defined in 2.9 with 
$S(\Fh^*)$ in 4.10.
We shall construct a standard $\BH$-module $E_{v,\r'}$ 
obtained from $H^{M^0(y_{\la})}_*(\CP_{y_{\la}}, \dot\CL^*)$ 
for a certain pair $(v, \r')$.
Let $y$ be the nilpotent element in $\Fg_{\BC}$ corresponding 
to $y_{\la} \in \Fg$.  
We choose $y^-, h_0 \in \Fg_{\BC}$ such that 
$\{ y, y^-, h_0\}$ forms  an $\Fs\Fl_2$-triple.
Put $h = h_0, r_0 = 1$.
Then $(h, r_0) \in \Fm(y)_{\BC}$ with $h$ semisimple.
We denote by $v$ an element in $H^*_{M^0(y_{\la})}$ corresponding to 
the $M^0(y)$-orbit of $(h,r_0)$.
Let $\r$ be the irreducible character of 
$A_G(y_{\la})$ as in 4.3.  Since 
$A_G(y_{\la}) \simeq \ol M(y_{\la}) \simeq \ol M(y)$, one can 
regard $\r$ as a character of $\ol M(y)$.  Let $\r^*$ be the 
dual representation of $\r$.
\par 
Under the notation in Remark 2.14 and Theorem 2.15, we note that 
\par\medskip\noindent
(4.12.1) \  
  Let 
$v$ be as above.  Then $E_{v,\r'}$ is a simple $\BH$-module, 
where $\r'$ is the restriction of $\r^*$ on $\ol M(y,v)$.
\par\medskip
It is enough to show that $(h, r_0)$ satisfies the
property in Theorem 2.15.   
By Dynkin-Kostant theory, $y$ is contained in the open dense orbit 
in $Y_{(h, r_0)} = \Fg_2$ (the graded space with respect to $h$) 
under the action of $Z_{G_{\BC}}(h)$. 
It remains to show that $\r'$ occurs in 
$H^{\{ e\}}_*(\CP_{y_{\la}}, \dot\CL^*)$.  But this is clear since
$H^{m}_c(\CP_{y_{\la}}, \dot\CL) = 
      H^{m}_c(\CP_{y_{\la}}, \dot\CL)_{\r}$. 
Thus (4.12.1) holds.
\para{4.13.}
The $\Fq$-structure $\vf_0: F^*\CL \isom \CL$ induces 
a linear isomorphism 
$\F$ on $H^*_c(\CP_{y_{\la}}, \dot\CL)$.  $\vf_0$ also induces a
linear isomorphism $\Psi$ on 
$H^{M^0(y_{\la})}_*(\CP_{y_{\la}},\dot\CL^*)$ satisfying the following 
property; by [L6, 7.2, (d)], there exists a $\Ql$-linear
isomorphism
\begin{equation*}
\tag{4.13.1}
\Ql\otimes_{H^*_{M^0(y_{\la})}}
      H_*^{M^0(y_{\la})}(\CP_{y_{\la}},\dot\CL^*) \to 
           H^{\{ e\}}_*(\CP_{y_{\la}},\dot\CL^*),
\end{equation*}
where $\Ql$ is regarded as an $H^*_{M^0(y_{\la})}$-module via
the canonical map $H^*_{M^0(y_{\la})} \to H^*_{\{ e\}} = \Ql$.
$F$ acts naturally on $H^*_{M^0(y_{\la})}$ and on $H^*_{\{ e\}}$, 
and the last map is $F$-equivariant with respect to the trivial action 
on $\Ql$. Thus $\Psi$ induces a linear map $\ol\Psi$ on the left hand
side of (4.13.1).   The $\Ql$-linear map in (4.13.1) is compatible 
with $\ol\Psi$ and the map $\F^*$ on 
$H^{\{ e\}}_*(\CP_{y_{\la}},\dot\CL^*) =
                 H^*_c(\CP_{y_{\la}},\dot\CL)^*$, 
where $\F^*$ is the transposed inverse of $\F$.
\par 
Note that $H^{m}_c(\CP_{y_{\la}},\dot\CL)$ is an irreducible 
$\CW$-module.  Since 
$\F w_0$ commutes with all the elements in $\CW$, we see that
$\F w_0$ acts on 
$H^{m}_c(\CP_{y_{\la}}, \dot\CL)$ as a scalar multiplication. 
Then we have the following lemma.
\begin{lem} 
Assume that $\F w_0$ acts on 
$H^{m}_c(\CP_{y_{\la}},\dot\CL)$ by a scalar 
multiplication by $\z$.  Then 
$\F w_0$ acts on $H^0_c(\CP_{y_{\la}},\dot\CL)$
by a scalar multiplication by $\z(-q\iv)^{m/2}$. 
\end{lem}
\begin{proof}
Let $v = (h, r_0)$ be as in (4.12.1).
Let $\g_v: H^*_{M^0(y_{\la})} \to \Ql$ be the algebra homomorphism 
corresponding to $v$ (cf. 2.13).  Since $M^0(y_{\la})$ is $F$-stable, 
$F$ acts naturally on $H^*_{M^0(y_{\la})}$ such that 
$\Psi(mx) = F(m)\Psi(x)$ for $m \in H^*_{M^0(y_{\la})}$ and 
$x \in H^{M^0(y_{\la})}_*(\CP_{y_{\la}},\dot\CL^*)$.  
Since $M_G^0(y_{\la})\simeq Z_G^0(y_{\la})\times \BG_m$,
we have 
$H_{M^0(y_{\la})}^* \simeq S(\Fh^*_1)^{W_1}\otimes \Ql[\Br]$, 
where $W_1$ is the Weyl group of a reductive group 
$Z^0_G(\f_{\la})$ and $\Fh_1^* = \Ql\otimes_{\BZ} X(T_1)$
with a maximally split maximal  torus $T_1$ of $Z^0_G(\f_{\la})$.
We note that
\par\medskip\noindent
(4.14.1) \ The maximal ideal $\Ker \g_v$ in $H^*_{M^0(y_{\la})}$
is generated by homogeneous polynomials.
\par\medskip
In fact, by the previous argument, we may replace 
$H^*_{M^0(y_{\la})}$ by $S(\Fm(y)^*_{\BC,r})^{M^0(y)}$, and 
$v$ by $(h, r_0) \in \Fm(y)_{\BC,r}$. 
It is enough to show that if a polynomial function $f$ on 
$\Fm(y)_{\BC,r}$ which is invariant under the action of $M^0(y)$ 
vanishes on $(h,r_0)$, then its homogeneous parts also vanish at
$(h,r_0)$.  But the $\BG_m$-action on $\Fm(y)_{\BC}$ implies that
$t : (h,r_0) \mapsto (t^{-2}h, t^{-2}r_0)$.  Since $f$ is invariant
under $M^0(y)$, we see that $f$ vanishes also on 
$(t^{-2}h,t^{-2}r_0)$ for any $t \in \BC^*$.  It follows that 
each homogeneous part of $f$ also vanishes at $(h,r_0)$ as asserted. 
\par
Next we note that 
\par\medskip\noindent
(4.14.2) \  The maximal ideal $\Ker \g_v$ is $F$-stable.
\par\medskip
Let $w_1$ be the longest element in $W_1$.  As in Lemma 4.11, 
$Fw_1$ acts on $S(\Fh_1^*)_i$ (the $i$-th homogeneous part) as a
scalar multiplication by $(-q)^i$.  Hence $F$ acts on 
$S(\Fh_1^*)^{W_1}_i$ by a scalar multiplication by $(-q)^i$. 
Also, $F$ acts on $\Ql[\Br]_i$ as a scalar multiplication by $q^i$.
Since $\Ker \g_v$ is a homogeneous ideal, $F$ stabilizes 
$\Ker\g_v$.  Hence (4.14.2) holds.
\par
Now $E_{v,\r'}$ is obtained as the quotient of 
$H^{M^0(y_{\la})}_*(\CP_{y_{\la}},\dot\CL^*)_{\r'}$ by the 
$\BH$-submodule
$I_v = \Ker \g_v\cdot H^{M^0(y_{\la})}_*(\CP_{y_{\la}},\dot\CL^*)_{\r'} $.
Since $\Ker\g_v$ is $F$-stable, we see that $I_v$ is $\Psi$-stable.  
Thus $\Psi$ induces a linear map
on $E_{v,\r'}$.
We consider the filtration $F^0 \subseteq F^1 \subseteq \cdots$ of
$H^{M^0(y_{\la})}_*(\CP_{y_{\la}},\dot\CL^*)$ as in 2.16.  Then
each $F^i$, as well as its $\r'$-isotypic part $F^i_{\r'}$,  
is $\Psi$-stable.  Then  $(F^i_{\r'})_v$ is also $\Psi$-stable since
it is the quotient of 
$F^i_{\r'}$ by $F^i_{\r'}\cap I_v$.
By (2.16.5) and by our assumption,  $\Psi w_0$ acts on 
the non-zero space $F_v^0 = (F^0_{\r'})_v $ as 
a scalar multiplication by $\z\iv$.
This implies that $\Psi w_0$ acts on 
$H^{M^0(y_{\la})}_0(\CP_{y_{\la}},\dot\CL^*)$ modulo $I_v$
by $\z\iv$.  On the other hand, since $E_{v,\r'}$ is a simple
$\BH$-module, $H^{M^0(y_{\la})}_*(\CP_{y_{\la}}, \dot\CL^*)_{\r'}$
is generated by $H^{M^0(y_{\la})}_0(\CP_{y_{\la}}, \dot\CL^*)_{\r'}$
mod $I_v$ as an $\BH$-module.  Since $\Br$ acts as a scalar
multiplication by $r_0$ on $E_{v,\r'}$, the action of $\BH$ on
$E_{v,\r'}$ is given by the action of $S(\Fh^*) = H^*_{T_0}$ 
and of $\CW$.  
Note that $\Psi w_0(\xi x) = Fw_0(\xi)\Psi w_0(x)$ for 
$\xi \in \BS, x \in H^{M^0(y_{\la})}_*(\CP_{y_{\la}},\dot\CL^*)$. 
The action of $\CW$ preserves the grading of
$H^{M^0(y_{\la})}_*(\CP_{y_{\la}},\dot\CL^*)$, and $\Psi w_0$ 
commutes with $\CW$.  It follows, by Lemma 4.11 that $\Psi w_0$
acts on $H^{M^0(y_{\la})}_{m}(\CP_{y_{\la}},\dot\CL^*)_{\r'}$
modulo $I_v$ as a scalar multiplication by $\z\iv(-q)^{m/2}$.
Let $f_m$ be the map $F^{m-1}_v \to F^m_v$ as in 2.16, which is
$\ol M(y_{\la}, v)$-equivariant.
Since 
$(F^{m}_v)_{\r'}/(\Im f_{m})_{\r'}$ 
is regarded as a natural quotient 
of $H^{M^0(y_{\la})}_{m}(\CP_{y_{\la}}, \dot\CL^*)_{\r'}$ 
modulo $I_v$, 
$\Psi w_0$ acts on $(F_v^{m})_{\r'}/(\Im f_{m})_{\r'}$ 
as $\z\iv(-q)^{m/2}$.  
Since $H^{\{ e\}}_{m}(\CP_{y_{\la}},\dot\CL^*)$ is isomorphic 
to $F^{m}_v/\Im f_{m}$ by (2.16.4),
we see that $\Psi w_0$ acts on 
$H^{\{ e\}}_{m}(\CP_{y_{\la}},\dot\CL^*)_{\r'}$ 
 by a scalar
multiplication by $\z\iv(-q)^{m/2}$, which coincides with the 
action of $\F^*w_0$ on it.
We claim that 
$H^{\{ e\}}_{m}(\CP_{y_{\la}},\dot\CL^*) = 
   H^{\{ e\}}_{m}(\CP_{y_{\la}},\dot\CL^*)_{\r'}$.
In fact, 
\begin{equation*}
H^{\{ e\}}_{m}(\CP_{y_{\la}}, \dot\CL^*) = 
H^0_c(\CP_{y_{\la}},\dot\CL)^* = \vG(\ZC, \CL)^* 
\end{equation*}
by Lemma 4.9.  $A_L(y_0)$ acts on $\vG(\ZC, \CL)$ by the character 
$\r_0$.  Since $\r$ is the pull back of $\r_0$ under the map 
$A_G(y_{\la}) \to A_L(y_0)$ (cf. 3.2), the action of 
$\ol M(y_{\la}) = A_G(y_{\la})$ on $H_c^0(\CP_{y_{\la}}, \dot\CL)$ 
is via $\r_0$.  Hence 
$H_c^0(\CP_{y_{\la}},\dot\CL) = H_c^0(\CP_{y_{\la}}, \dot\CL)_{\r}$  
and the claim follows.
\par
Thus $\F w_0$ acts on 
$H^0_c(\CP_{y_{\la}},\dot\CL) = 
        H^{\{ e\}}_{m}(\CP_{y_{\la}}, \dot\CL^*)^*$
by a scalar multiplication by $\z(-q\iv)^{m/2}$ as asserted.
\end{proof}
\para{4.15.}
We are now ready to prove Theorem 4.4.  First we note that 
$\CW$ acts trivially on $H^0_c(\CP_y, \dot\CL)$ for any 
$y = y_{\nu}$ such that all the parts of 
$\nu$ are divisible by $d$.
In fact, if $y_{\nu}$ is regular nilpotent, $a_0+r = 0$ by 
(1.3.1) since $\ZC$ is also a regular nilpotent class in $L$.
It follows, by the generalized Springer correspondence (see 3.2), 
that $H^0_c(\CP_{y_{\nu}}, \dot\CL)$ is the irreducible $\CW$-module    
corresponding to the unit representation.   Thus by Lemma 4.9,
$H^0_c(\CP_y, \dot\CL)$ is also a trivial $\CW$-module for any
$y$.
\par
Now assume that $\F w_0$ acts on 
$H^{m}_c(\CP_{y_{\la}},\dot\CL)$ by a scalar multiplication 
by $\z$.  Then by Lemma 4.14, $\F w_0$ acts on $H^0_c(\CP_{y_{\la}},\dot\CL)$
by a scalar multiplication by $\z(-q\iv)^{m/2}$.
Since the map 
$H^0_c(\CP_{y_1}, \dot\CL) \to H^0_c(\CP_{y_{\la}},\dot\CL)$ 
is $\F w_0$-equivariant isomorphism by Lemma 4.9 (and Proposition
4.8), we see that
$\F w_0$ acts on $H^0_c(\CP_{y_1}, \dot\CL)$ by a scalar
multiplication by $\z(-q\iv)^{m/2}$. 
Since $w_0$ acts trivially on it, we see that  
\par\medskip\noindent
(4.15.1) \ 
$\F$ acts on 
$H^0_c(\CP_{y_1}, \dot\CL)$ by $\z(-q\iv)^{m/2}$.
\par\medskip 
On the other hand, by a similar argument as in the proof of
Lemma 4.9, the natural map
\begin{equation*}
H^0_c(\CP_{y_0}, \dot\CL) \simeq \vG(\ZC, \CL) \to \CL_{y_0}
\end{equation*}
gives an isomorphism.  This isomorphism is 
compatible with the action of $\F$ and of $\vf_0$.  It follows that
$\F$ acts on $H_c^0(\CP_{y_0}, \dot\CL)$ as an identity map. 
Since $y_1$ is in the $G$-orbit of $y_0$, 
$H^0_c(\CP_{y_0},\dot\CL) \simeq H^0_c(\CP_{y_1}, \dot\CL)$.  
As discussed in the proof of Lemma 4.14, $A_G(y_0)$ acts on 
$H_c^0(\CP_{y_0},\dot\CL)$ via $\r_0$.  We also note that 
$A_L(y_0) \simeq A_G(y_0)$. 
Since $y_1$ is $G^F$-conjugate to
$y_{c_1}$, an element twisted by $c_1 \in A_G(y_0)$ from $y_0$,
we see that  
\par\medskip\noindent
(4.15.2) \ $\F$ acts on $H^0_c(\CP_{y_1},\dot\CL)$ by a scalar
multiplication by $\r_0(c_1) = \e_{\la}$.
\par\medskip
Comparing (4.15.1) and (4.15.2), we see that  
$\z = \e_{\la}(-q)^{m/2}$.  This proves the theorem.
\para{4.16.}
In order to apply Theorem 4.4, we need to know $c_1 \in A_G(y_0)$
such that $y_1 = (y_0)_{c_1}$.  For a given $y_0$, we shall choose 
a specific $y_1$ and $y_{\la}$, and determnie $c_1$ explicitly. 
Put $\la' = (\la_1', \dots, \la_r')$ with $\la_i' = \la_i/d$.
Hence $\la'$ is a partition of $t$.  Let $\{e_j^{(i)}\}$ be the basis 
of $V_0$ as in 4.3.  Put $d' = [d/2]$.  Let us define a subspace 
$W_0$ of $V_0$ and define $y_0$ by 
\begin{equation*}
W_0 = \begin{cases}
          \lp e^{(i)}_{d'+1} \mid 1 \le i \le t \rp
                &\quad\text{ if $d$ is odd},  \\
          \{ 0\} &\quad\text{ if $d$ is even}.
      \end{cases}
\end{equation*}
Also we define subspaces $W_1, W_2$ of $V_0$ by
\begin{align*}
W_1 &= \lp e^{(i)}_j \mid 1 \le j \le d', 1 \le i \le t\rp, \\
W_2 &= \lp e^{(i)}_j \mid d-d'+1 \le j \le d, 1 \le i \le t\rp.
\end{align*}
Clearly we have $V_0 = W_1\oplus W_0\oplus W_2$.  
We define a new basis $\{ h^{(i)}_j \mid 1\le j \le d,1 \le i \le t\}$
of $V_0$ satisfying the following conditions.
\par\medskip
\begin{enumerate}
\item
$h_j^{(i)} = e_j^{(i)}$ if $e_j^{(i)} \in W_1$. 
\item
The set $\{ h^{(i)}_j \mid d-d'+1 \le j \le d, 1 \le i \le t\}$  
coincides with the set of the basis $\{ e^{(i)}_j\}$ of $W_2$.
\item
Let $z$ be the number of $i$ such that $\la_i'$ is odd. Then 
\begin{align*}
\lp h_j^{(2i-1)}, h_{d-j+1}^{(2i)}\rp &= 1 
    \quad\text{ for } 1 \le i \le (t-z)/2, 1 \le j \le d, \\
\lp h_j^{(i)}, h_{d-j+1}^{(i)} \rp &= 1 
    \quad\text{ for }(t-z)/2 +1 \le i \le t, 1 \le j \le d',\\
\lp h^{(i)}_{d'+1}, h^{(i)}_{d'+1}\rp &= \pm 1 
    \quad\text{ for } (t-z)/2 +1 \le i \le t \text{ if $d$: odd}, \\
\lp h^{(i)}_j, h^{(i')}_{j'}\rp &= 0.
    \quad\text{ for all other cases.}  
\end{align*}
\item
$\{ h^{(i)}_{d'+1}\}$ gives a basis of $W_0$ in the case where $d$ is odd.
\end{enumerate}
The conditions (i) $\sim$ (iv) determines $\{ h^{(i)}_j\}$ uniquely 
except the vectors contained in $W_0$.  We note that one can 
choose the basis $\{ h^{(i)}_j\}$ of $W_0$ so that the transition
matrix between $\{e^{(i)}_j\}$ and $\{h^{(i)}_j\}$ has the determinant
$\pm 1$ (see the construction of $f_j^{(i)}$ in 4.2).
\par
We define a nilpotent transformation $y'_1 \in \Fg^{F}$ as in 4.2 
replacing $f_j^{(i)}$ by $h_j^{(i)}$.   Then it is easy to construct
$y'_{\la} \in \Fg^F$ such that $y_1'$ is obtained from $y_{\la}'$ by 
a similar procedure as $y_1$ is obtained from $y_{\la}$.  Clearly, 
$y_1'$ (resp. $y_{\la}'$) is conjugate to $y_1$ (resp. $y_{\la}$). 
The argument in the proof of Theorem 4.4 works well for 
$y_1', y'_{\la}$.  Thus we consider $y_1'$ and $y_{\la}'$, and write
them as $y_1$, $y_{\la}$.  We shall describe $c_1 \in A_G(y_0)$ 
such that $y_1 = (y_0)_{c_1}$.  It follows from the construction that 
there  
exists $g \in \wt G = GL_n$ such that $\Ad(g) y_0 = y_1$, where 
$g$ stabilizes the subspaces $W_0, W_1, W_2$.  
More precisely, $g$ acts trivially on $W_1$, and gives a permutation 
matrix with respect to the basis $\{ e^{(i)}_j\}$ up to $\pm 1$ on 
$W_2$.  Thus by our choice of the basis $\{ h^{(i)}_j\}$, we have
$\det g = \pm 1$.  Let us take $\a \in \ol\BF_q$ such that 
$\a^n = \pm 1$ (if $\det g = 1$, we choose $\a = 1$). 
We note that $\Ker y_0 \subset W_2$, and that $g$ stabilizes 
$\Ker y_0$.  We denote by $g_0$ the restriction of $g$ on $\Ker y_0$.
Then we have
\begin{lem}  
Let the notations be as before.   Then  we have
$y_1 = (y_0)_{c_1}$, where $c_1 \in  A_G(y_0)$ is given, 
under the identification 
$A_G(y_0) \simeq \{ x \in \ol\BF_q^* \mid x^d = 1\}$, by  
\begin{equation*}
c_1 = \a^{t(1-q)}\det g_0.
\end{equation*}
\end{lem}
\begin{proof}
Let $\f_0: \Fs\Fl_2 \to \Fg$ be as in 4.10, and we consider 
the group $Z_G(\f_0)$.  Then 
$A_G(y_0) \simeq Z_G(\f_0)/Z^0_G(\f_0)$.
We have $Z_G(\f_0) \simeq \{ x \in GL_t \mid \det x^d = 1\}$, where
the element $g_1 \in Z_G(\f_0)$ corresponding to $x$ is given as
follows; $g_1$ acts on the subspace $V_j$ of $V_0$ spanned by 
$\{ e^{(i)}_j \mid 1 \le i \le t \}$, for a fixed $j$, as  $x \in GL_t$. 
Now if we can find $g_1 \in G$ such that $Ad(g_1)y_0 = y_1$, then 
$g_1\iv F(g_1) \in A_G(y_0)$, and it leaves $\Ker y_0$ invariant.
Moreover, the determinant of the restriction of $g_1\iv F(g_1)$ gives
rise to the corresponding element in 
$A_G(y_0) \simeq \{ x \in \ol\BF_q^* \mid x^d = 1\}$. 
\par 
Now in our situation, if we put $g_1 = \a\iv g$, we have 
$g_1 \in G$
and $\Ad(g_1)y_0 = y_1$.  Then $g_1\iv F(g_1) = \a^{1-q} g\iv F(g)$.
On the other hand, since $F(g) = w_0({}^tg\iv) w_0\iv$, 
$F(g)$ also stabilizes the subspaces $W_0, W_1, W_2$.  Moreover, 
it acts on $W_2$ trivially, and on $W_1$ as a permutaiton of the
basis $\{ e^{(i)}_j\}$ up to sign.  It follows that $g\iv F(g)$
acts on the space $\Ker y_0$ as $g_0\iv$.  Thus $g_1\iv F(g_1)$ 
acts on $\Ker y_0$ as a map $\a^{1-q}g_0\iv$, and we have 
$\det (\a^{1-q}g_0\iv) = \a^{t(1-q)}\det g_0$ as asserted.
\end{proof}

\par\medskip
\par

\end{document}